\documentclass[11pt]{article}
\usepackage{amsmath,amsfonts,amsthm,amssymb}
\usepackage{cite}
\usepackage{indentfirst}
\usepackage{needspace}
\usepackage[colorlinks=true, linkcolor=blue, citecolor=blue, urlcolor=blue]{hyperref}
\hypersetup{pdftitle={A Breadth-First Search Characterization of 2-Threshold Graphs},
	pdfauthor={Qinze Zhu and Yingzhi Tian},pdflang={en}}
\newtheoremstyle{citedplain}{\topsep}{\topsep}{\itshape}{}{\bfseries}{.}{.5em}
{\thmname{#1}\thmnumber{ #2}\thmnote{ \normalfont#3}}
\theoremstyle{citedplain}
\newtheorem{theorem}{Theorem}[section]
\newtheorem{lemma}[theorem]{Lemma}
\newtheorem{corollary}[theorem]{Corollary}
\newtheorem{claim}{Claim}
\theoremstyle{definition}
\newtheorem{definition}{Definition}[section]
\theoremstyle{citedplain}
\makeatletter
\newcommand{\Rmnum}[1]{\expandafter\@slowromancap\romannumeral #1@}
\makeatother

\allowdisplaybreaks[3]
\begin{document}
	\title{A Breadth-First Search Characterization of $2$-Threshold Graphs\footnote{The research is supported by Natural Science Foundation of Xinjiang Uygur Autonomous Region (2025D01E02) and National Natural Science Foundation of China (12261086).}}
	\author{Qinze Zhu, Yingzhi Tian\footnote{Corresponding author. E-mail: qinzezhu@stu.xju.edu.cn (Q. Zhu); tianyzhxj@163.com (Y. Tian).} \\
		{\small College of Mathematics and System Sciences, Xinjiang
			University, Urumqi, Xinjiang 830046, PR China}}
	\date{}
	\maketitle

	\begin{sloppypar}
		\clubpenalty=10000
		\widowpenalty=10000
		\noindent{\bf Abstract }
		A graph is $2$-threshold if its vertices admit real ranks and two
		thresholds $\theta_1<\theta_2$ such that distinct vertices are adjacent
		exactly when their rank sum belongs to $[\theta_1,\theta_2)$.
		Jamison and Sprague proved that bipartite $2$-threshold graphs are
		precisely the bipartite permutation graphs.
		We characterize connected non-bipartite $2$-threshold graphs by
		conditions on breadth-first search (BFS) levels and prove sufficiency
		by an explicit rank construction.
		We also characterize disconnected $2$-threshold graphs using
		the known component criterion.
		
		\noindent{\bf Keywords:} Multithreshold graphs; $2$-threshold graphs; Breadth-first search

		\section{Introduction}
		All graphs in this paper are finite, simple, and undirected.  The class of
		threshold graphs was introduced by Chv\'{a}tal and Hammer \cite{Chvatal1} in
		connection with integer linear programming.  A graph $G$ is a threshold graph if
		there exist a real-valued rank function $r:V(G)\to\mathbb R$ and a threshold
		$\theta$ such that, for any two distinct vertices $u$ and $v$,
		$uv\in E(G)\Longleftrightarrow r(u)+r(v)\ge\theta$.
		Threshold graphs admit many equivalent descriptions. A nonempty graph is
		a threshold graph if and only if it can be obtained from one vertex by
		repeatedly adding an isolated vertex or a vertex adjacent to every existing
		vertex.  Equivalently, its vertices can be
		partitioned into a clique and a stable set so that the neighborhoods of the
		stable-set vertices inside the clique are linearly ordered by inclusion.  These
		characterizations and their algorithmic consequences are discussed in
		\cite{Golumbic}.
		
		Threshold graphs also arise in the study of synchronization primitives
		for parallel processing, where Henderson and Zalcstein independently
		obtained the same class \cite{Henderson}.
		Ordman studied their use in resource allocation \cite{Ordman85}
		and the memory requirements for synchronization \cite{Ordman89}.
		Applications to graph labelling for open shop scheduling \cite{EckerZaks}
		and cyclic workforce scheduling \cite{Koop} are described in the monograph
		\cite{MahadevPeled}.
		Threshold graphs and related classes also occur in the study of Guttman
		scales and biorders in mathematical psychology
		\cite{Guttman,Doignon,MahadevPeled}.
		Diaconis, Holmes and Janson studied random threshold graphs and their
		limits \cite{Diaconis}.
		
		Jamison and Sprague \cite{Jamison} introduced multithreshold graphs.
		In their model, each vertex has a real rank, and two distinct vertices
		are adjacent exactly when the number of thresholds at or below their
		rank sum is odd. For a positive integer $k$, a $k$-threshold representation
		of a graph $G$ consists of thresholds $\theta_1<\theta_2<\cdots<\theta_k$
		and a rank function $r:V(G)\to\mathbb R$ such that, for distinct $u,v$,
		$uv\in E(G)\Longleftrightarrow
		r(u)+r(v)\in[\theta_{2i-1},\theta_{2i})$
		for some $i\in\{1,\ldots,\lceil \frac{k}{2}\rceil\}$,
		where $\theta_{k+1}=+\infty$ when $k$ is odd.
		A graph admitting such a representation is called $k$-threshold.
		The case $k=1$ gives the classical threshold graphs. For $k=2$,
		the adjacency condition is $\theta_1\le r(u)+r(v)<\theta_2$.
		The threshold number of $G$, denoted by
		$\Theta(G)$, is the least $k$ for which $G$ has a $k$-threshold representation.
		Restricting a representation gives
		$\Theta(H)\le \Theta(G)$ for every induced subgraph $H$ of $G$.
		
		Jamison and Sprague \cite{Jamison} proved that every $n$-vertex graph with $n\ge2$ satisfies
		$\Theta(G)\le \binom n2$ and constructed graphs with threshold number linear in
		$n$.  They also asked for the threshold number of complete multipartite graphs.
		Puleo proved the lower bound
		$\Theta(nK_3)\ge n^{\frac{1}{3}}$ \cite{Puleo}.
		Chen and Hao \cite{Chen} obtained exact results for complete multipartite graphs
		with sufficiently large parts and for their complements.
		Kittipassorn and Sumalroj \cite{Kittipassorn} determined the threshold numbers of
		$K_{3,3,\ldots,3}$, $K_{4,4,\ldots,4}$ and their complements $nK_3$ and $nK_4$.
		Here $K_{n_1,\ldots,n_h}$ denotes the graph with stable parts of sizes
		$n_1,\ldots,n_h$ and all edges between distinct parts, and $nK_q$
		denotes the disjoint union of $n$ complete graphs on $q$ vertices.
		
		$2$-threshold graphs are precisely the pairwise compatibility graphs
		that admit a representation on a star \cite{Kobayashi,Monti}. Xiao and Nagamochi \cite{Xiao}
		characterized this class by gap-free vertex orderings and gave an
		$O(n^3m)$-time recognition algorithm, where $n$ and $m$ are the
		numbers of vertices and edges.
		
		A permutation graph is the intersection graph of line segments joining
		two parallel lines. Jamison and Sprague \cite{Jamison1} proved that every
		$2$-threshold graph is a permutation graph and that the bipartite
		$2$-threshold graphs are precisely the bipartite permutation graphs.
		Kobayashi, Okamoto, Otachi and Uno \cite{Kobayashi} gave a linear-time
		recognition algorithm for $2$-threshold graphs, which they call
		\emph{double-threshold graphs}. For a non-bipartite graph $G$, let $M$
		be a maximum clique minimizing the sum of the degrees of its vertices in
		$G$. They proved that $G$ is $2$-threshold if and only if the auxiliary
		graph they construct from $G$ and $M$ is a bipartite permutation graph.
		
		Enright, Stewart and Tardos \cite{Enright} proved that every connected
		permutation graph has a BFS root with the following property. For
		consecutive levels $L_i,L_{i+1}$ and vertices $x,y\in L_i$, either $y$
		is adjacent to every neighbor of $x$ in $L_{i+1}$, or $x$ is adjacent
		to every neighbor of $y$ in $L_{i+1}$. Jamison and Sprague \cite{Jamison1}
		used BFS levels to construct representations of bipartite permutation graphs.
		
		We characterize connected non-bipartite $2$-threshold graphs using the
		BFS levels $L_0,\ldots,L_d$ of the original graph. For a suitable root,
		edges within a level occur only in $L_{d-1}$ or $L_d$. Theorem~$\ref{the:necessary}$ specifies
		which pairs $x,y\in L_i$ must satisfy
		$N_{L_{i-1}}(y)\subseteq N_{L_{i-1}}(x)$. From a root and level
		decompositions satisfying these conditions, we construct a $2$-threshold
		representation.
		
		Section 2 gives the notation and preliminary lemmas.
		Section~$\ref{sec:characterization}$ derives the necessary conditions from
		BFS levels, and Section~$\ref{sec:construction}$ constructs a rank representation.
		Section 5 treats bipartite and disconnected graphs.
		Section~$\ref{sec:conclusion}$ gives concluding remarks.
		
		\section{Preliminaries}
		For a graph $G$,
		let $V(G)$ and $E(G)$ denote its vertex and edge sets. Write
		$uv$ for the edge with distinct endpoints $u,v$. The open
		neighborhood of $v$ is $N_G(v)=\{w\in V(G):vw\in E(G)\}$,
		and its closed neighborhood is $N_G[v]=N_G(v)\cup\{v\}$.
		A vertex $v$ is isolated if $N_G(v)=\emptyset$. For subsets $X,Y\subseteq V(G)$,
		put $N_G(X)=\left(\bigcup\limits_{v\in X}N_G(v)\right)\setminus X$.
		Write $N_Y(X)=N_G(X)\cap Y$.
		For a vertex $v$, write $N_Y(v)=N_G(v)\cap Y$; when $v\in Y$,
		write $N_Y[v]=N_Y(v)\cup\{v\}$. These definitions allow $X$ and
		$Y$ to intersect. The induced subgraph on $X$ is $G[X]$, and
		$G-X=G[V(G)\setminus X]$.
		
		A set $X$ is a clique if every pair of distinct vertices in $X$
		is adjacent, and it is a stable set if $E(G[X])=\emptyset$.
		A maximum clique is a clique of largest cardinality; its size
		is the clique number $\omega(G)$.
		
		A family $\mathcal C$ is nested if, for any $X,Y\in\mathcal C$,
		either $X\subseteq Y$ or $Y\subseteq X$.
		Its distinct members can then be listed as a strict inclusion chain.
		
		In an ordering $x_1,\ldots,x_m$, we say that $x_j$ precedes
		$x_k$ if $j<k$. For $0\le k\le m$, the set
		$\{x_1,\ldots,x_k\}$ is an initial segment, and
		$\{x_{k+1},\ldots,x_m\}$ is a terminal segment.
		Both the empty set and the full set are allowed in either case.
		An interval in this order
		is either empty or has the form $\{x_j,x_{j+1},\ldots,x_k\}$
		for some $1\le j\le k\le m$.
		
		We write $X\mathbin{\dot\cup}Y$ for the union of disjoint sets $X,Y$.
		A graph is bipartite if its vertex set has a partition
		$V(G)=X\mathbin{\dot\cup}Y$ into two stable sets.
		
		For a connected graph, $d_G(u,v)$ is the length of a shortest
		$u$--$v$ path and
		$\operatorname{ecc}_G(u)=\max\limits_{v\in V(G)}d_G(u,v)$.
		
		We order the vertices within each BFS level by rank and use strict
		representations to exclude ties among ranks and rank sums,
		as well as equalities involving thresholds.
		\begin{definition}[Strict representation]\label{def:strict}
			Following Jamison and Sprague \cite{Jamison}, a $k$-threshold
			representation is strict if
			
			(1) the vertex ranks are pairwise distinct;
			
			(2) the sums $r(u)+r(v)$, indexed by unordered pairs of distinct
			vertices, are pairwise distinct;
			
			(3) no vertex rank equals a threshold;
			
			(4) no pairwise rank sum equals a threshold;
			
			(5) no vertex rank equals a pairwise rank sum.
		\end{definition}
		
		\begin{lemma}[{\cite{Jamison}}]\label{lem:strict}
			For any $k$-threshold representation of a graph, there exists a
			strict $k$-threshold representation of the same graph with the
			same thresholds.
		\end{lemma}
		
		In a strict representation, the rank order is the increasing order
		of the vertex ranks.
		
		\begin{lemma}[{\cite{Jamison1}}]\label{lem:neighborhood}
			Let $(r,(\theta_1,\theta_2))$ be a strict $2$-threshold
			representation of $G$. For each $p\in V(G)$, the set $N_G(p)$
			is an interval in the rank order on $V(G)\setminus\{p\}$:
			if $u,v,w\in V(G)\setminus\{p\}$ satisfy $u,w\in N_G(p)$ and
			$r(u)<r(v)<r(w)$, then $v\in N_G(p)$.
		\end{lemma}
		
		\begin{lemma}\label{lem:extreme-adjacent}
			Let $G$ be a $2$-threshold graph with no isolated vertices, and let
			$(r,(\theta_1,\theta_2))$ be a $2$-threshold representation of $G$.  If $u$ and
			$v$ are distinct vertices of minimum and maximum rank, respectively, then $uv\in E(G)$.
		\end{lemma}
		\noindent{\bf Proof.} Suppose, to the contrary, that $uv\notin E(G)$.  Since $G$ has no isolated
		vertices, choose $u'\in N_G(u)$ and $v'\in N_G(v)$.
		Since $r(u')\le r(v)$ and $r(u)\le r(v')$, the edges $uu'$ and $vv'$ give
		$\theta_1\le r(u)+r(u')\le r(u)+r(v)\le r(v')+r(v)<\theta_2$.
		Thus $uv\in E(G)$, a contradiction. $\Box$
		
		For a $2$-threshold graph with no isolated vertices, take a strict
		representation and let $v$ have maximum rank.
		Lemma~$\ref{lem:extreme-adjacent}$ places the minimum-rank vertex in
		$N_G(v)$, and Lemma~$\ref{lem:neighborhood}$ makes $N_G(v)$ an interval
		of $V(G)\setminus\{v\}$. It is therefore an initial segment of that
		rank order. To compare neighborhoods at subsequent steps, we use
		the following two inclusions obtained from a common neighbor of
		maximum or minimum rank.
		
		\begin{lemma}\label{lem:rank-nesting}
			Let $G$ be a $2$-threshold graph with representation
			$(r,(\theta_1,\theta_2))$. Let $X\subseteq V(G)$ and let $u,v,w\in X$ be
			distinct vertices such that $r(u)\le r(v)$ and
			$w\in N_X(u)\cap N_X(v)$.
			
			(i) If $r(w)\ge \max\{r(\xi):\xi\in (N_X(u)\cup N_X(v))\setminus\{u,v\}\}$, then
			$N_X(u)\subseteq N_X[v]$.
			
			(ii) If $r(w)\le \min\{r(\xi):\xi\in (N_X(u)\cup N_X(v))\setminus\{u,v\}\}$, then
			$N_X(v)\subseteq N_X[u]$.
		\end{lemma}
		\noindent{\bf Proof.} For (i), suppose that there is a vertex
		$u'\in N_X(u)\setminus N_X[v]$.
		Since $r(u')\le r(w)$ and $r(u)\le r(v)$, the edges $uu'$ and $vw$ give
		$\theta_1\le r(u)+r(u')\le r(v)+r(u')\le r(v)+r(w)<\theta_2$.
		$vu'\in E(G)$, contrary to $u'\notin N_X[v]$, so
		$N_X(u)\subseteq N_X[v]$.
		
		For (ii), suppose that $v'\in N_X(v)\setminus N_X[u]$.
		Since $r(w)\le r(v')$, the edges $uw$ and $vv'$ give
		$\theta_1\le r(u)+r(w)\le r(u)+r(v')\le r(v)+r(v')<\theta_2$.
		$uv'\in E(G)$, contrary to $v'\notin N_X[u]$, so
		$N_X(v)\subseteq N_X[u]$. $\Box$

		\section{BFS layers and structural properties}\label{sec:characterization}
		We derive necessary conditions by taking a vertex of maximum or minimum
		rank as the BFS root.
		
		\subsection{Levels and next-level neighborhoods}\label{subsec:bfs-partitions}
		\begin{definition}[BFS levels]\label{def:bfs}
			Let $G$ be connected, and let $u\in V(G)$.
			The breadth-first search (BFS) levels are
			$L_i=\{v\in V(G):d_G(u,v)=i\}$ for $0\le i\le d$, where
			$d=\operatorname{ecc}_G(u)$. We call $u$ the root.
			For a vertex in $L_i$, the levels $L_{i-1}$ and $L_{i+1}$,
			when they exist, are its preceding and next levels.
			An edge with both endpoints in $L_i$ is a same-level edge.
		\end{definition}
		
		Then $L_0=\{u\}$, and each of $L_0,\ldots,L_d$ is nonempty.
		For $i\ge1$, a shortest path from $u$ to a vertex of $L_i$ gives
		that vertex a neighbor in $L_{i-1}$.
		Every edge joins vertices in the same level or in consecutive levels.
		Consequently, every clique is contained in one level or in the union
		of two consecutive levels.
		If all levels are stable, the unions of the even and odd levels
		form a bipartition.
		
		\begin{theorem}\label{the:level-nesting}
			Let $G$ be a connected $2$-threshold graph with at least two vertices,
			let $(r,(\theta_1,\theta_2))$ be a strict representation, and take a vertex $u$
			of maximum rank as the BFS root. Let $L_0,\ldots,L_d$ be the resulting BFS levels and
			put $U_i=\bigcup\limits_{j=i}^dL_j$.
			
			(i) For $1\le i\le d$, the set $L_i$ is an initial segment of $U_i$
			in increasing rank order when $i$ is odd, and a terminal segment
			when $i$ is even.
			
			(ii) For $0\le i<d$ and $x,y\in L_i$ with $r(x)<r(y)$,
			$N_{L_{i+1}}(x)\subseteq N_{L_{i+1}}(y)$ when $i$ is odd, and
			$N_{L_{i+1}}(y)\subseteq N_{L_{i+1}}(x)$ when $i$ is even.
			
			(iii) For $1\le i\le d$ and $x,y\in L_i$ with $r(x)<r(y)$,
			$N_{L_{i-1}}(y)\subseteq N_{L_{i-1}}(x)$ when $i$ is odd, and
			$N_{L_{i-1}}(x)\subseteq N_{L_{i-1}}(y)$ when $i$ is even.
		\end{theorem}
		
		\noindent\textbf{Proof.} Order the vertices $u_1,\ldots,u_n$ by
		$r(u_1)<\cdots<r(u_n)$, so $u=u_n$.
		For (i), we proceed by induction on the level index $i$.
		Lemma~$\ref{lem:extreme-adjacent}$ gives $u_1u_n\in E(G)$.
		Then Lemma~$\ref{lem:neighborhood}$ yields
		$L_1=N_G(u_n)=\{u_1,\ldots,u_{\ell_1}\}$, where $\ell_1=|L_1|$.
		This proves (i) for $i=1$.
		
		Assume $d\ge2$ and put $X_1=\{u\}\cup L_1\cup L_2$.
		The root $u$ is adjacent to every vertex of $L_1$ and has maximum rank in $X_1$.
		Since $u_{\ell_1}$ has maximum rank in $L_1$, Lemma~$\ref{lem:rank-nesting}$\,(i) gives
		$N_{X_1}(x_1)\subseteq N_{X_1}[u_{\ell_1}]$ for $x_1\in L_1\setminus\{u_{\ell_1}\}$,
		and hence $N_{L_2}(x_1)\subseteq N_{L_2}(u_{\ell_1})$.
		Every vertex of $L_2$ has a neighbor in $L_1$, so
		$N_{L_2}(u_{\ell_1})=\bigcup\limits_{x\in L_1}N_{L_2}(x)=L_2$.
		
		Take $u_a\in L_2$ and $u_b\in U_2$ with $a<b$.
		Both $u_a$ and $u$ are neighbors of $u_{\ell_1}$, and
		$r(u_a)<r(u_b)<r(u)$. Lemma~$\ref{lem:neighborhood}$ gives
		$u_{\ell_1}u_b\in E(G)$. Since $u_{\ell_1}\in L_1$ and $u_b\in U_2$, the BFS
		definition forces $u_b\in L_2$. Hence $L_2$ is a terminal
		segment of $U_2$.
		
		Assume $d\ge3$, and let $u_{k_1}$ have minimum rank in $L_2$.
		Since $N_{L_2}(u_{\ell_1})=L_2$ and $\max\limits_{x\in L_1}r(x)<\min\limits_{y\in U_2}r(y)$,
		the vertex $u_{\ell_1}$ has minimum rank in
		$X_2=\{u_{\ell_1}\}\cup L_2\cup L_3$.
		Lemma~$\ref{lem:rank-nesting}$\,(ii) gives
		$N_{X_2}(x_2)\subseteq N_{X_2}[u_{k_1}]$ for $x_2\in L_2\setminus\{u_{k_1}\}$,
		and hence $N_{L_3}(x_2)\subseteq N_{L_3}(u_{k_1})$.
		Since every vertex of $L_3$ has a neighbor in $L_2$,
		$N_{L_3}(u_{k_1})=\bigcup\limits_{x\in L_2}N_{L_3}(x)=L_3$.
		
		Take $u_{a_2}\in L_3$ and $u_{b_2}\in U_3$ with $b_2<a_2$.
		Then $r(u_{\ell_1})<r(u_{b_2})<r(u_{a_2})$.
		Both $u_{\ell_1}$ and $u_{a_2}$ are neighbors of $u_{k_1}$, so
		Lemma~$\ref{lem:neighborhood}$ gives $u_{k_1}u_{b_2}\in E(G)$.
		Since $u_{b_2}\in U_3$, it follows that $u_{b_2}\in L_3$.
		Consequently, $L_3$ is an initial segment of $U_3$.
		
		For $d\ge4$, start the induction at $i=3$, with
		$v_1=u_{\ell_1}$ and $v_2=u_{k_1}$ from the base cases.
		For $3\le i<d$, assume (i) for $1\le j\le i$ and that vertices
		$v_j\in L_j$, $1\le j<i$, have been chosen with maximum rank in $L_j$
		when $j$ is odd and minimum rank when $j$ is even, satisfying
		$N_{L_{j+1}}(v_j)=L_{j+1}$.
		Put $X_i=\{v_{i-1}\}\cup L_i\cup L_{i+1}$.
		
		If $i$ is odd, then $L_{i-1}$ is a terminal segment of $U_{i-1}$,
		so $\max\limits_{x\in U_i}r(x)<\min\limits_{y\in L_{i-1}}r(y)$,
		and $v_{i-1}$ has maximum rank in $X_i$.
		Choose $v_i\in L_i$ with maximum rank.
		For $x_3\in L_i\setminus\{v_i\}$, Lemma~$\ref{lem:rank-nesting}$\,(i) gives
		$N_{X_i}(x_3)\subseteq N_{X_i}[v_i]$, and hence
		$N_{L_{i+1}}(x_3)\subseteq N_{L_{i+1}}(v_i)$.
		Every vertex of $L_{i+1}$ has a neighbor in $L_i$, so
		$N_{L_{i+1}}(v_i)=\bigcup\limits_{x\in L_i}N_{L_{i+1}}(x)=L_{i+1}$.
		Take $u_{a_3}\in L_{i+1}$ and $u_{b_3}\in U_{i+1}$ with $a_3<b_3$.
		By (i) for $L_{i-1}$, $r(u_{a_3})<r(u_{b_3})<r(v_{i-1})$.
		Since $u_{a_3},v_{i-1}\in N_G(v_i)$, Lemma~$\ref{lem:neighborhood}$ gives
		$v_iu_{b_3}\in E(G)$. As $u_{b_3}\in U_{i+1}$ and $v_i\in L_i$,
		this implies $u_{b_3}\in L_{i+1}$.
		Thus $L_{i+1}$ is a terminal segment of $U_{i+1}$.
		
		For even $i$, $L_{i-1}$ is an initial segment of $U_{i-1}$,
		so $\max\limits_{x\in L_{i-1}}r(x)<\min\limits_{y\in U_i}r(y)$,
		and $v_{i-1}$ has minimum rank in $X_i$.
		Choose $v_i\in L_i$ with minimum rank.
		Analogously to the odd case, using Lemma~$\ref{lem:rank-nesting}$\,(ii) in place of (i),
		$N_{L_{i+1}}(v_i)=L_{i+1}$.
		Take $u_{a_4}\in L_{i+1}$ and $u_{b_4}\in U_{i+1}$ with $b_4<a_4$.
		By (i) for $L_{i-1}$, $r(v_{i-1})<r(u_{b_4})<r(u_{a_4})$.
		The same application of Lemma~$\ref{lem:neighborhood}$ gives
		$v_iu_{b_4}\in E(G)$ and hence $u_{b_4}\in L_{i+1}$.
		Thus $L_{i+1}$ is an initial segment of $U_{i+1}$.
		In both cases, $v_i$ has the required properties, completing the induction for (i).
		
		For (ii), the case $i=0$ is vacuous. Let $1\le i<d$ and take
		$x_4,y_4\in L_i$ with $r(x_4)<r(y_4)$.
		In $X_i$ constructed in the proof of (i), these vertices have a common neighbor
		of maximum rank for odd $i$ and minimum rank for even $i$.
		For odd $i$, the same application of Lemma~$\ref{lem:rank-nesting}$\,(i) gives
		$N_{X_i}(x_4)\subseteq N_{X_i}[y_4]$, and hence
		$N_{L_{i+1}}(x_4)\subseteq N_{L_{i+1}}(y_4)$.
		For even $i$, analogously using Lemma~$\ref{lem:rank-nesting}$\,(ii),
		$N_{L_{i+1}}(y_4)\subseteq N_{L_{i+1}}(x_4)$.
		
		It remains to prove (iii).
		For $i=1$, the preceding-level neighborhoods are all $\{u\}$,
		which gives (iii) in this case.
		
		For $2\le i\le d$, take $x_5,y_5\in L_i$ with $r(x_5)<r(y_5)$.
		By (ii) for $i-2$, choose $u_{k_4}\in L_{i-2}$ whose neighborhood in $L_{i-1}$
		is largest under inclusion.
		Since every vertex of $L_{i-1}$ has a neighbor in $L_{i-2}$,
		$N_{L_{i-1}}(u_{k_4})=\bigcup\limits_{v\in L_{i-2}}N_{L_{i-1}}(v)=L_{i-1}$.
		If $i$ is odd, (i) for $L_{i-2}$ gives
		$r(u_{k_4})<r(x_5)<r(y_5)$.
		Take $u_c\in N_{L_{i-1}}(y_5)$.
		Then $y_5,u_{k_4}\in N_G(u_c)$; Lemma~$\ref{lem:neighborhood}$ gives
		$u_cx_5\in E(G)$, so $N_{L_{i-1}}(y_5)\subseteq N_{L_{i-1}}(x_5)$.
		
		For even $i$, (i) for $L_{i-2}$, or the maximum rank of the root $u$
		when $i=2$, gives $r(x_5)<r(y_5)<r(u_{k_4})$.
		Analogously to the odd case, interchanging $x_5$ and $y_5$ gives
		$N_{L_{i-1}}(x_5)\subseteq N_{L_{i-1}}(y_5)$.
		$\Box$
		
		If instead the BFS root $u$ has minimum rank, put
		$\widetilde r(v)=-r(v)$ for every $v\in V(G)$ and
		$(\widetilde\theta_1,\widetilde\theta_2)=(-\theta_2,-\theta_1)$.
		Then $\widetilde\theta_1<\widetilde\theta_2$, and strictness gives,
		for distinct $x,y$,
		\[
		xy\in E(G)\Longleftrightarrow
		\theta_1<r(x)+r(y)<\theta_2
		\Longleftrightarrow
		\widetilde\theta_1<\widetilde r(x)+\widetilde r(y)<\widetilde\theta_2.
		\]
		The negated ranks and thresholds still satisfy all conditions of a strict representation.
		Thus $(\widetilde r,(\widetilde\theta_1,\widetilde\theta_2))$ is a strict
		representation of the same graph, in which $u$ has maximum rank.
		For this fixed root, the BFS levels are unchanged.
		Theorem~$\ref{the:level-nesting}$\,(i) therefore makes $L_i$ a terminal
		segment of $U_i$ for odd $i$ and an initial segment for even $i$
		in the original increasing rank order.
		For $r(x)<r(y)$, we have $\widetilde r(y)<\widetilde r(x)$;
		applying (ii) and (iii) to $(y,x)$ reverses their neighborhood
		inclusions for each parity of $i$.
		
		For the rest of this section, when the BFS root has maximum rank,
		we order $L_i$ by increasing rank for odd $i$ and by decreasing rank
		for even $i$. When the root has minimum rank, we use decreasing rank
		for odd $i$ and increasing rank for even $i$.
		The latter orders are exactly the former orders for
		$(-r,(-\theta_2,-\theta_1))$.
		
		For BFS levels $L_0,\ldots,L_d$,
		$\bigcup\limits_{v\in L_i}N_{L_{i+1}}(v)=L_{i+1}$ for $0\le i<d$;
		if these neighborhoods are nested, one of them equals $L_{i+1}$.
		Together with Theorem~$\ref{the:level-nesting}$\,(ii)--(iii), the level orders
		fixed above, and the rank-negation argument, this gives the following corollary.
		
		\Needspace{4\baselineskip}
		\begin{corollary}\label{cor:next-nesting}
			Let $G$ be a connected graph with at least two vertices and
			BFS levels $L_0,\ldots,L_d$ rooted at $u$.
			
			(i) For $0\le i<d$, if $\{N_{L_{i+1}}(v):v\in L_i\}$ is nested,
			then there exists a vertex $v\in L_i$ such that $N_{L_{i+1}}(v)=L_{i+1}$.
			
			(ii) If $G$ has a strict $2$-threshold representation
			$(r,(\theta_1,\theta_2))$ and $u$ has maximum or minimum rank,
			then, for distinct $x,y\in L_i$ with $x$ preceding $y$,
			$N_{L_{i+1}}(x)\subseteq N_{L_{i+1}}(y)$ when $i<d$, and
			$N_{L_{i-1}}(y)\subseteq N_{L_{i-1}}(x)$ when $i\ge1$.
			In particular, $\{N_{L_{i+1}}(v):v\in L_i\}$ is nested
			for every $0\le i<d$.
		\end{corollary}
		
		In each level, we group vertices with the same next-level neighborhood.
		
		\begin{definition}[Next-level neighborhood classes]\label{def:next-classes}
			Let $G$ be a connected graph with at least two vertices and
			BFS levels $L_0,\ldots,L_d$ rooted at $u$ such that $\{N_{L_{i+1}}(v):v\in L_i\}$ is nested for every
			$0\le i<d$. For each such $i$, choose representatives
			$u^{(i)}_1,\ldots,u^{(i)}_{s_i}\in L_i$, one for each distinct
			nonempty next-level neighborhood, indexed so that
			$N_{L_{i+1}}(u^{(i)}_1)\subsetneqq\cdots\subsetneqq
			N_{L_{i+1}}(u^{(i)}_{s_i})$.
			Put $I_i=\{1,\ldots,s_i\}$ and define
			$D_{i,t}=\{v\in L_i:N_{L_{i+1}}(v)
			=N_{L_{i+1}}(u^{(i)}_t)\}$ for $t\in I_i$, and
			$D_{i,0}=\{v\in L_i:N_{L_{i+1}}(v)=\emptyset\}$.
		\end{definition}
		
		By Corollary~$\ref{cor:next-nesting}$\,(i) and the neighborhood ordering in
		Definition~$\ref{def:next-classes}$, $s_i\ge1$ and
		$N_{L_{i+1}}(u^{(i)}_{s_i})=L_{i+1}$.
		At the root, $s_0=1$ and $u^{(0)}_1=u$.
		
		The sets $D_{i,t}$, $0\le t\le s_i$, satisfy
		$L_i=D_{i,0}\mathbin{\dot\cup}D_{i,1}\mathbin{\dot\cup}\cdots\mathbin{\dot\cup}D_{i,s_i}$,
		and only $D_{i,0}$ may be empty.
		
		In the level orders fixed after Theorem~$\ref{the:level-nesting}$,
		each class $D_{i,t}$ is an interval, and the nonempty classes occur
		in increasing order of $t$.
		If $x$ precedes $y$ and $y$ precedes $z$, with $x,z\in D_{i,t}$,
		then Corollary~$\ref{cor:next-nesting}$\,(ii) gives
		$N_{L_{i+1}}(x)\subseteq N_{L_{i+1}}(y)
		\subseteq N_{L_{i+1}}(z)=N_{L_{i+1}}(x)$, so $y\in D_{i,t}$.
		The same inclusion orders distinct classes by increasing index.
		
		\subsection{Same-level edges and admissible decompositions}
		We use the following decompositions of the BFS levels.
		
		\Needspace{3\baselineskip}
		\begin{definition}[Admissible decomposition]\label{def:admissible}
			A decomposition $L_i=A_i\mathbin{\dot\cup}Q_i$ of a BFS level
			is admissible if it has the following properties:
			\begin{enumerate}
				\item[(i)] $A_i$ is a nonempty stable set, and $Q_i$ is a clique,
				possibly empty;
				\item[(ii)] the family $\{N_{L_i}(x):x\in A_i\}$ is nested;
				\item[(iii)] $N_{L_i}(A_i)=Q_i$.
			\end{enumerate}
		\end{definition}
		
		If $G[L_i]$ is edgeless, Definition~$\ref{def:admissible}$\,(iii) forces $Q_i=\emptyset$ and
		$A_i=L_i$, which is admissible.
		Conversely, $Q_i=\emptyset$ implies that $L_i=A_i$ is stable.
		Hence $Q_i=\emptyset$ if and only if
		$G[L_i]$ is edgeless, and the admissible decomposition in this
		case is unique. In particular, $A_0=\{u\}$ and $Q_0=\emptyset$.
		
		\Needspace{5\baselineskip}
		\begin{lemma}\label{lem:same-level}
			Let $G$ be a connected $2$-threshold graph with at least two vertices,
			let $(r,(\theta_1,\theta_2))$ be a strict representation, and take a vertex $u$
			of maximum or minimum rank as the BFS root.
			
			(i) If $x\in L_i$ with $0\le i<d$ has a same-level neighbor,
			then $i=d-1$ and $N_{L_d}(x)=L_d$.
			In particular, every $L_i$ with $0\le i\le d-2$ is stable.
			
			(ii) For $0\le i\le d$ and distinct $x,y\in L_i$ with $x$ preceding $y$,
			$N_{L_i}(x)\setminus\{y\}\subseteq N_{L_i}(y)\setminus\{x\}$.
			
			(iii) Write the order on each $L_i$ as $x_1,\ldots,x_m$, and let
			$k$ be the largest index for which $\{x_1,\ldots,x_k\}$ is stable.
			Then $A_i=\{x_1,\ldots,x_k\}$ and $Q_i=L_i\setminus A_i$
			give an \mbox{admissible decomposition}.
			
			(iv) For the decompositions in (iii),
			$Q_{d-1}\cup N_{A_{d-1}}(Q_{d-1})\subseteq D_{d-1,s_{d-1}}$.
		\end{lemma}
		
		\noindent\textbf{Proof.} Negating the ranks as described after
		Theorem~$\ref{the:level-nesting}$ preserves the level orders and neighborhoods,
		and hence the sets $A_i,Q_i,D_{i,t}$ constructed from them.
		It suffices to assume that $u$ has maximum rank, with odd levels
		ordered by increasing rank and even levels by decreasing rank.
		
		For (i), suppose that $1\le i<d$ and that $v_1\in L_i$ has a
		same-level neighbor $v_2$. Put $U_{i+1}=\bigcup\limits_{j=i+1}^dL_j$,
		choose $w_1\in N_{L_{i-1}}(v_1)$, and take any $w_2\in U_{i+1}$.
		
		For odd $i$, Theorem~$\ref{the:level-nesting}$\,(i) gives
		$r(v_2)<r(w_2)<r(w_1)$, where $w_1=u$ when $i=1$.
		Since $v_1v_2,v_1w_1\in E(G)$, strictness gives
		$\theta_1<r(v_1)+r(v_2)<r(v_1)+r(w_2)<r(v_1)+r(w_1)<\theta_2$.
		Thus $v_1w_2\in E(G)$ for every $w_2\in U_{i+1}$.
		Since $v_1\in L_i$, it follows that
		$i+1\le d_G(u,w_2)\le d_G(u,v_1)+1=i+1$.
		Hence $U_{i+1}=L_{i+1}$, $i=d-1$, and $N_{L_d}(v_1)=L_d$.
		
		For even $i$, Theorem~$\ref{the:level-nesting}$\,(i) gives
		$r(w_1)<r(w_2)<r(v_2)$, so
		$\theta_1<r(v_1)+r(w_1)<r(v_1)+r(w_2)<r(v_1)+r(v_2)<\theta_2$.
		Thus $v_1w_2\in E(G)$ for every $w_2\in U_{i+1}$.
		The same distance inequalities give $U_{i+1}=L_{i+1}$,
		$i=d-1$, and $N_{L_d}(v_1)=L_d$.
		
		For $1\le i\le d-2$, an edge in $G[L_i]$ would force $i=d-1$,
		a contradiction. Hence each such $L_i$ is stable, as is $L_0=\{u\}$.
		
		For (ii), the case $i=0$ is vacuous.
		Fix $1\le i\le d$.
		Definition~$\ref{def:next-classes}$ and Corollary~$\ref{cor:next-nesting}$\,(i) give
		$N_{L_i}(u^{(i-1)}_{s_{i-1}})=L_i$.
		Put $Y_i=L_i\cup\{u^{(i-1)}_{s_{i-1}}\}$ and take distinct $x,y\in L_i$
		with $x$ preceding $y$ in this level order.
		If $i$ is odd, Theorem~$\ref{the:level-nesting}$\,(i) gives
		$r(v)<r(u^{(i-1)}_{s_{i-1}})$ for every $v\in L_i$.
		For $i=1$, this follows from $u^{(0)}_1=u$ and the maximum rank of the root.
		Thus $r(x)<r(y)$ and $u^{(i-1)}_{s_{i-1}}$ is a common
		neighbor of maximum rank in $Y_i$. Lemma~$\ref{lem:rank-nesting}$\,(i) gives
		$N_{Y_i}(x)\subseteq N_{Y_i}[y]$, and hence
		$N_{L_i}(x)\setminus\{y\}\subseteq N_{L_i}(y)\setminus\{x\}$.
		
		For even $i$, Theorem~$\ref{the:level-nesting}$\,(i) gives
		$r(u^{(i-1)}_{s_{i-1}})<r(v)$ for every $v\in L_i$.
		Analogously to the odd case, using Lemma~$\ref{lem:rank-nesting}$\,(ii)
		with $(y,x)$ in place of $(x,y)$ gives
		$N_{L_i}(x)\setminus\{y\}\subseteq N_{L_i}(y)\setminus\{x\}$.
		
		For (iii), fix $0\le i\le d$ and write its level order as
		$x_1,\ldots,x_m$. The index $k$ exists and satisfies $k\ge1$
		because $\{x_1\}$ is stable. If $k=m$, then $G[L_i]$ is edgeless,
		and $A_i=L_i$, $Q_i=\emptyset$ is admissible by
		Definition~$\ref{def:admissible}$. This includes every $i\le d-2$
		by (i).
		
		Suppose $k<m$. By the choice of $k$, $A_i$ is nonempty and stable.
		Since $k$ is maximal, the vertex $x_{k+1}$ has a neighbor
		$x_{a_1}$ with $a_1\le k$.
		If $a_1=k$, then $x_kx_{k+1}\in E(G)$.
		If $a_1<k$, the inclusion in (ii) for $x_{a_1},x_k$, applied to
		$x_{k+1}$, gives the same edge.
		For $b_1>k+1$, apply (ii) to $x_{k+1},x_{b_1}$ with the neighbor
		$x_k$; it follows that $x_kx_{b_1}\in E(G)$.
		Thus $N_{L_i}(x_k)=Q_i$.
		For any $k<a_2<b_2$, apply (ii) to $x_k,x_{a_2}$ with the neighbor
		$x_{b_2}$ to obtain $x_{a_2}x_{b_2}\in E(G)$. Hence $Q_i$ is a clique.
		Since $A_i$ is stable, (ii) reduces to
		$N_{L_i}(x_{a_3})\subseteq N_{L_i}(x_{b_3})$ for $a_3<b_3\le k$.
		These neighborhoods are nested and contained in $Q_i$, and
		$N_{L_i}(x_k)=Q_i$ gives $N_{L_i}(A_i)=Q_i$, so the decomposition is admissible.
		
		For (iv), since $A_{d-1}$ is stable, its vertices with same-level
		neighbors are exactly $N_{A_{d-1}}(Q_{d-1})$.
		For each $q\in Q_{d-1}$, Definition~$\ref{def:admissible}$\,(iii) gives
		$v\in A_{d-1}$ with $vq\in E(G)$.
		Thus the vertices with same-level neighbors are exactly
		$Q_{d-1}\cup N_{A_{d-1}}(Q_{d-1})$.
		By (i), each has next-level neighborhood $L_d$, so
		$Q_{d-1}\cup N_{A_{d-1}}(Q_{d-1})\subseteq D_{d-1,s_{d-1}}$. $\Box$
		
		\subsection{Stable classes and a maximum clique}\label{subsec:stable-classes}
		For $0\le i\le d-2$, the level $L_i=A_i$ is stable, and we keep
		its next-level classes $D_{i,t}$ unchanged.
		In $L_{d-1}$, we further group only the stable vertices of the largest class
		$D_{d-1,s_{d-1}}$ by their same-level neighborhoods.
		In $L_d$, we use only same-level neighborhoods.
		The notation $S_{i,t}$ is used only in these last two levels.
		
		\Needspace{4\baselineskip}
		\begin{definition}[Stable neighborhood classes]\label{def:stable-classes}
			Let $G$ be a connected graph with at least two vertices and
			BFS levels $L_0,\ldots,L_d$ rooted at $u$ and nested next-level neighborhoods.
			Use the classes in Definition~$\ref{def:next-classes}$ and admissible
			decompositions $L_i=A_i\mathbin{\dot\cup}Q_i$ with $Q_i=\emptyset$
			for $0\le i\le d-2$ and
			$Q_{d-1}\cup N_{A_{d-1}}(Q_{d-1})\subseteq D_{d-1,s_{d-1}}$.
			
			(i) For $i=d-1$, partition the nonempty set
			$A_{d-1}\cap D_{d-1,s_{d-1}}$ according to equality of the
			same-level neighborhoods. Choose representatives
			$u^{(d-1)}_{s_{d-1}},\ldots,u^{(d-1)}_{m_{d-1}}$ for the distinct
			neighborhoods, including the empty one when it occurs, and order them
			so that
			\[
			N_{L_{d-1}}(u^{(d-1)}_{s_{d-1}})\subsetneqq\cdots\subsetneqq
			N_{L_{d-1}}(u^{(d-1)}_{m_{d-1}}).
			\]
			Choosing a new representative $u^{(d-1)}_{s_{d-1}}$ does not change
			the classes $D_{d-1,t}$.
			For $s_{d-1}\le t\le m_{d-1}$, define
			\[
			S_{d-1,t}=\{v\in A_{d-1}\cap D_{d-1,s_{d-1}}:
			N_{L_{d-1}}(v)=N_{L_{d-1}}(u^{(d-1)}_t)\}.
			\]
			For $0\le t<s_{d-1}$, put $S_{d-1,t}=D_{d-1,t}$.
			
			(ii) For $i=d$, if $Q_d\ne\emptyset$, choose representatives
			$u^{(d)}_1,\ldots,u^{(d)}_{m_d}\in A_d$ for the distinct nonempty
			same-level neighborhoods, in increasing inclusion order. Put
			$S_{d,t}=\{v\in A_d:N_{L_d}(v)=N_{L_d}(u^{(d)}_t)\}$
			for $1\le t\le m_d$, and
			$S_{d,0}=\{v\in A_d:N_{L_d}(v)=\emptyset\}$.
			If $Q_d=\emptyset$, set $m_d=0$, $S_{d,0}=L_d$, and choose
			$u^{(d)}_0\in L_d$.
		\end{definition}
		
		We call $D_{i,t}$ for $0\le i\le d-2$ and $S_{i,t}$ for
		$i\in\{d-1,d\}$ the stable classes.
		The nonempty classes partition $A_i$, with $A_i=L_i$ for $i\le d-2$.
		
		\begin{corollary}\label{cor:terminal-neighborhoods}
			Let $G$ be a connected graph with at least two vertices and
			BFS levels $L_0,\ldots,L_d$ rooted at $u$ and nested next-level neighborhoods.
			Suppose each $L_i=A_i\mathbin{\dot\cup}Q_i$ is admissible,
			$Q_i=\emptyset$ for $0\le i\le d-2$, and
			$Q_{d-1}\cup N_{A_{d-1}}(Q_{d-1})\subseteq D_{d-1,s_{d-1}}$.
			For each $i\in\{d-1,d\}$, the following hold.
			
			(i) $N_{L_i}(u^{(i)}_{m_i})=Q_i$, including when $Q_i=\emptyset$.
			
			(ii) $Q_i\cup\{u^{(i)}_{m_i}\}$ is a maximum clique of $G[L_i]$,
			and $\omega(G[L_i])=|Q_i|+1$.
		\end{corollary}
		
		\noindent\textbf{Proof.} Choose the representatives as in
		Definition~$\ref{def:stable-classes}$.
		For (i), if $Q_i=\emptyset$, then $L_i=A_i$ is stable and
		$N_{L_i}(u^{(i)}_{m_i})=\emptyset=Q_i$.
		Otherwise, by Definition~$\ref{def:admissible}$\,(ii)--(iii),
		the same-level neighborhoods on $A_i$ are nested and satisfy
		$\bigcup\limits_{v\in A_i}N_{L_i}(v)=Q_i$.
		Choose $x\in A_i$ with $N_{L_i}(x)=Q_i$.
		For $i=d-1$, $Q_i\ne\emptyset$ gives
		$x\in N_{A_{d-1}}(Q_{d-1})\subseteq D_{d-1,s_{d-1}}$.
		Since $u^{(i)}_{m_i}$ has the largest same-level neighborhood
		among the representatives,
		$Q_i=N_{L_i}(x)\subseteq N_{L_i}(u^{(i)}_{m_i})\subseteq Q_i$.
		
		For (ii), by (i) and Definition~$\ref{def:admissible}$\,(i),
		$Q_i\cup\{u^{(i)}_{m_i}\}$ is a clique, while every clique
		meets the stable set $A_i$ in at most one vertex.
		Thus $\omega(G[L_i])=|Q_i|+1$. $\Box$
		
		Lemma~3.6 of \cite{Kobayashi} implies that a non-bipartite $2$-threshold
		graph has a maximum clique whose deletion leaves a bipartite graph.
		We obtain such a clique by combining maximum cliques of $G[L_{d-1}]$ and
		$G[L_d]$.
		
		\begin{lemma}\label{lem:terminal-clique}
			Let $G$ be a connected non-bipartite graph with BFS levels
			$L_0,\ldots,L_d$ rooted at $u$ and nested next-level neighborhoods.
			Suppose each $L_i=A_i\mathbin{\dot\cup}Q_i$ is admissible,
			$Q_i=\emptyset$ for $0\le i\le d-2$, and
			$Q_{d-1}\cup N_{A_{d-1}}(Q_{d-1})\subseteq D_{d-1,s_{d-1}}$.
			For $i\in\{d-1,d\}$, let $K_i=Q_i\cup\{u^{(i)}_{m_i}\}$.
			Then $Q_{d-1}\cup Q_d\ne\emptyset$, $K=K_{d-1}\cup K_d$ is a maximum
			clique of $G$, and $G-K$ is bipartite. Moreover,
			$\omega(G)=\omega(G[L_{d-1}])+\omega(G[L_d])=|Q_{d-1}|+|Q_d|+2$.
		\end{lemma}
		
		\noindent\textbf{Proof.} Set $Q=Q_{d-1}\cup Q_d$. Since $Q_i=\emptyset$ for
		$0\le i\le d-2$, deleting $Q$ leaves the stable set $A_i$ in each level $L_i$.
		Let $X$ and $Y$ be the sets of remaining vertices in the even and odd levels,
		respectively. Every edge of $G-Q$ joins $X$ to $Y$, so $G-Q$ is bipartite. In particular, $Q\ne\emptyset$ because $G$ is non-bipartite.
		Since $Q\subseteq K$, $G-K$ is an induced subgraph of $G-Q$ and is also bipartite.
		
		By Corollary~$\ref{cor:terminal-neighborhoods}$\,(ii), $K_i$ is a maximum
		clique of $G[L_i]$ for $i\in\{d-1,d\}$. The hypothesis gives
		$Q_{d-1}\subseteq D_{d-1,s_{d-1}}$. By Definition~$\ref{def:stable-classes}$\,(i),
		$u^{(d-1)}_{m_{d-1}}\in A_{d-1}\cap D_{d-1,s_{d-1}}$. Therefore
		$K_{d-1}=Q_{d-1}\cup\{u^{(d-1)}_{m_{d-1}}\}\subseteq D_{d-1,s_{d-1}}$.
		By Definition~$\ref{def:next-classes}$ and Corollary~$\ref{cor:next-nesting}$\,(i),
		each vertex of $K_{d-1}$ is adjacent to every vertex in $L_d$. Hence $K$ is a clique and
		$|K|\le\omega(G)$.
		
		To prove $\omega(G)\le |K|$, let $K'$ be a maximum clique of $G$.
		We have $|K'\cap Q|\le |Q|$ since $K'\cap Q\subseteq Q$.
		The vertices of $K'$ are pairwise adjacent and $X$ and $Y$ are stable,
		so $|K'\cap X|\le 1$ and $|K'\cap Y|\le 1$.
		Then $|K'|=|K'\cap Q|+|K'\cap X|+|K'\cap Y|\le |Q|+2$.
		Since $K=Q\mathbin{\dot\cup}\{u^{(d-1)}_{m_{d-1}},u^{(d)}_{m_d}\}$,
		$\omega(G)=|K'|\le |Q|+2=|K|$. Thus $\omega(G)=|K|=|K_{d-1}|+|K_d|
		=\omega(G[L_{d-1}])+\omega(G[L_d])=|Q_{d-1}|+|Q_d|+2$.~$\Box$
		
		\subsection{Neighborhood differences and inclusions}\label{subsec:compatibility}
		We first take successive differences of the next-level neighborhood
		chains, and then of the same-level chains in the last two levels.
		
		\begin{definition}[Preceding-level difference classes]\label{def:t-differences}
			For $1\le i\le d$, let $u^{(i-1)}_1,\ldots,u^{(i-1)}_{s_{i-1}}\in L_{i-1}$
			be the representatives of the distinct nonempty next-level
			neighborhoods in Definition~$\ref{def:next-classes}$,
			ordered by increasing inclusion. Put
			$T^{(i-1)}_1=N_{L_i}(u^{(i-1)}_1)$ and
			$T^{(i-1)}_j=N_{L_i}(u^{(i-1)}_j)\setminus N_{L_i}(u^{(i-1)}_{j-1})$
			for $2\le j\le s_{i-1}$.
		\end{definition}
		
		\Needspace{3\baselineskip}
		\begin{corollary}\label{cor:t-adjacency}
			For $1\le i\le d$, the following hold.
			
			(i) The sets $T^{(i-1)}_j$, $1\le j\le s_{i-1}$, are nonempty
			and partition $L_i$.
			
			(ii) For $v\in T^{(i-1)}_j$,
			$N_{L_{i-1}}(v)=\bigcup\limits_{t=j}^{s_{i-1}}D_{i-1,t}$.
			Consequently, for $x\in D_{i-1,t}$, $0\le t\le s_{i-1}$,
			we have $xv\in E(G)\Longleftrightarrow j\le t$.
			
			(iii) For $v_1\in T_j^{(i-1)}$ and $v_2\in T_{j'}^{(i-1)}$,
			$N_{L_{i-1}}(v_2)\subseteq N_{L_{i-1}}(v_1)\Longleftrightarrow j\le j'$.
			
			\Needspace{3\baselineskip}
			(iv) For the stable classes in Definition~$\ref{def:stable-classes}$,
			$D_{d-1,s_{d-1}}=Q_{d-1}\cup
			\bigcup\limits_{t=s_{d-1}}^{m_{d-1}}S_{d-1,t}$, and every vertex of this
			union is adjacent to all of $L_d$.
			For $0\le t<s_{d-1}$, $x\in S_{d-1,t}$, and $v\in T^{(d-1)}_j$,
			we have $xv\in E(G)\Longleftrightarrow j\le t$.
		\end{corollary}
		
		\noindent\textbf{Proof.} For (i), Corollary~$\ref{cor:next-nesting}$\,(i)
		and the neighborhood chain in Definition~$\ref{def:next-classes}$
		give nonempty successive differences with disjoint union $L_i$.
		
		For (ii), fix $v\in T^{(i-1)}_j$. For $1\le t\le s_{i-1}$, nesting gives
		$v\in N_{L_i}(u^{(i-1)}_t)\Longleftrightarrow j\le t$,
		and every vertex of $D_{i-1,t}$ has the same next-level neighborhood
		as $u^{(i-1)}_t$. Thus $xv\in E(G)\Longleftrightarrow j\le t$
		for $x\in D_{i-1,t}$. The same equivalence holds for $t=0$,
		since $D_{i-1,0}$ has no neighbors in $L_i$ and $j\ge1$.
		Taking the union over $t\ge j$ gives the stated formula for $N_{L_{i-1}}(v)$.
		
		For (iii), if $j\le j'$, (ii) gives
		$N_{L_{i-1}}(v_2)=\bigcup\limits_{t=j'}^{s_{i-1}}D_{i-1,t}
		\subseteq\bigcup\limits_{t=j}^{s_{i-1}}D_{i-1,t}=N_{L_{i-1}}(v_1)$.
		If $j>j'$, choose $z\in D_{i-1,j'}$, which is nonempty by
		Definition~$\ref{def:next-classes}$. By (ii),
		$z\in N_{L_{i-1}}(v_2)\setminus N_{L_{i-1}}(v_1)$, so
		$N_{L_{i-1}}(v_2)\not\subseteq N_{L_{i-1}}(v_1)$.
		
		For (iv), Definition~$\ref{def:stable-classes}$\,(i) partitions
		$A_{d-1}\cap D_{d-1,s_{d-1}}$ into the sets $S_{d-1,t}$ with
		$s_{d-1}\le t\le m_{d-1}$.
		Since $Q_{d-1}\subseteq D_{d-1,s_{d-1}}$,
		$D_{d-1,s_{d-1}}=Q_{d-1}\cup\bigcup\limits_{t=s_{d-1}}^{m_{d-1}}S_{d-1,t}$.
		For every $v\in D_{d-1,s_{d-1}}$, Definition~$\ref{def:next-classes}$
		and Corollary~$\ref{cor:next-nesting}$\,(i) give $N_{L_d}(v)=L_d$.
		For $t<s_{d-1}$, the equality $S_{d-1,t}=D_{d-1,t}$ and (ii)
		give the remaining adjacency assertion. $\Box$
		
		\begin{samepage}
			\begin{definition}[Clique difference classes]\label{def:r-differences}
				Use the representatives $u^{(i)}_t\in A_i$ of the distinct
				same-level neighborhoods from Definition~$\ref{def:stable-classes}$,
				ordered by increasing inclusion.
				Put $R^{(d-1)}_{s_{d-1}}=N_{L_{d-1}}(u^{(d-1)}_{s_{d-1}})$
				and, when $m_d\ge1$, $R^{(d)}_1=N_{L_d}(u^{(d)}_1)$.
				For $i=d-1$ and $s_{d-1}<t\le m_{d-1}$,
				or for $i=d$ and $2\le t\le m_d$, put
				$R^{(i)}_t=N_{L_i}(u^{(i)}_t)\setminus N_{L_i}(u^{(i)}_{t-1})$.
			\end{definition}
		\end{samepage}
		
		\Needspace{4\baselineskip}
		\begin{corollary}\label{cor:r-adjacency}
			For $i\in\{d-1,d\}$, the following hold.
			
			(i) The nonempty sets $R^{(i)}_t$ partition $Q_i$.
			Only $R^{(d-1)}_{s_{d-1}}$ can be empty.
			
			(ii) For $w\in R^{(i)}_t$ and $x\in S_{i,p}$ with $0\le p\le m_i$,
			$xw\in E(G)\Longleftrightarrow t\le p$, and
			$N_{A_i}(w)=\bigcup\limits_{p=t}^{m_i}S_{i,p}$.
		\end{corollary}
		
		\Needspace{3\baselineskip}
		\noindent\textbf{Proof.} For (i) with $i=d-1$, Corollary~$\ref{cor:terminal-neighborhoods}$\,(i)
		gives $N_{L_{d-1}}(u^{(d-1)}_{m_{d-1}})=Q_{d-1}$.
		Thus $Q_{d-1}=R^{(d-1)}_{s_{d-1}}\mathbin{\dot\cup}\cdots
		\mathbin{\dot\cup}R^{(d-1)}_{m_{d-1}}$.
		The strict inclusions in Definition~$\ref{def:stable-classes}$\,(i)
		make $R_t^{(d-1)}$ nonempty for $s_{d-1}<t\le m_{d-1}$.
		If $Q_{d-1}=\emptyset$, Definition~$\ref{def:stable-classes}$\,(i)
		gives $m_{d-1}=s_{d-1}$ and
		$R^{(d-1)}_{s_{d-1}}=\emptyset$.
		
		For $i=d$, if $Q_d=\emptyset$, Definition~$\ref{def:stable-classes}$\,(ii)
		gives $m_d=0$, so no set $R_t^{(d)}$ is defined.
		Otherwise, Definition~$\ref{def:stable-classes}$\,(ii) gives
		$R^{(d)}_1=N_{L_d}(u^{(d)}_1)\ne\emptyset$.
		Analogously to $i=d-1$,
		$Q_d=R^{(d)}_1\mathbin{\dot\cup}\cdots\mathbin{\dot\cup}R^{(d)}_{m_d}$,
		and every $R_t^{(d)}$ is nonempty.
		
		For (ii), first suppose that $s_{d-1}\le p\le m_{d-1}$ when $i=d-1$,
		or $1\le p\le m_d$ when $i=d$.
		By Definition~$\ref{def:stable-classes}$,
		$N_{L_i}(x)=N_{L_i}(u^{(i)}_p)$.
		Since $w\in R_t^{(i)}$, Definition~$\ref{def:r-differences}$ gives
		$xw\in E(G)\Longleftrightarrow w\in N_{L_i}(u^{(i)}_p)
		\Longleftrightarrow t\le p$.
		
		If $i=d-1$ and $p<s_{d-1}$, then $x\in D_{d-1,p}$, whereas
		$N_{A_{d-1}}(Q_{d-1})\subseteq D_{d-1,s_{d-1}}$.
		Thus $xw\notin E(G)$ and $p<s_{d-1}\le t$.
		
		For $i=d$ and $p=0$, Definition~$\ref{def:stable-classes}$\,(ii) gives
		$xw\notin E(G)$, and $p<t$, so the equivalence holds.
		Since the classes $S_{i,p}$ partition $A_i$,
		$N_{A_i}(w)=\bigcup\limits_{p=t}^{m_i}S_{i,p}$. $\Box$
		
		\Needspace{4\baselineskip}
		We next compare neighborhoods in the preceding level $L_{i-1}$.
		
		\begin{lemma}\label{lem:clique-intervals}
			Let $G$ be a connected $2$-threshold graph with at least two vertices,
			let $(r,(\theta_1,\theta_2))$ be a strict representation, and take a vertex $u$
			of maximum or minimum rank as the BFS root. Use the admissible decompositions $L_i=A_i\mathbin{\dot\cup}Q_i$
			constructed in Lemma~$\ref{lem:same-level}$\,(iii).
			For $i\in\{d-1,d\}$ with $i\ge1$, the following hold.
			
			(i) If $x_1\in A_i$ and $y_1\in Q_i$, then
			$N_{L_{i-1}}(y_1)\subseteq N_{L_{i-1}}(x_1)$.
			
			(ii) If $x_2\in R^{(i)}_{t'}$ and $y_2\in R^{(i)}_t$, with
			$s_{d-1}\le t<t'\le m_{d-1}$ when $i=d-1$ and
			$1\le t<t'\le m_d$ when $i=d$, then
			$N_{L_{i-1}}(y_2)\subseteq N_{L_{i-1}}(x_2)$.
		\end{lemma}
		
		\noindent\textbf{Proof.} For (i), the choice of $A_i,Q_i$ in
		Lemma~$\ref{lem:same-level}$\,(iii) puts $x_1$ before $y_1$.
		Corollary~$\ref{cor:next-nesting}$\,(ii) gives the required inclusion.
		
		For (ii), Corollary~$\ref{cor:r-adjacency}$\,(ii) gives
		$N_{A_i}(x_2)=\bigcup\limits_{p=t'}^{m_i}S_{i,p}
		\subsetneqq\bigcup\limits_{p=t}^{m_i}S_{i,p}=N_{A_i}(y_2)$.
		The inclusion is strict because
		$\emptyset\ne S_{i,t}\subseteq N_{A_i}(y_2)\setminus N_{A_i}(x_2)$.
		If $y_2$ preceded $x_2$, Lemma~$\ref{lem:same-level}$\,(ii) would give
		$N_{A_i}(y_2)\subseteq N_{A_i}(x_2)$, a contradiction.
		Thus $x_2$ precedes $y_2$, and Corollary~$\ref{cor:next-nesting}$\,(ii) gives
		$N_{L_{i-1}}(y_2)\subseteq N_{L_{i-1}}(x_2)$.
		$\Box$
		
		Corollary~$\ref{cor:next-nesting}$ gives nested neighborhoods between consecutive levels.
		Lemma~$\ref{lem:same-level}$ gives admissible decompositions of the levels, and
		Lemma~$\ref{lem:clique-intervals}$ compares neighborhoods in the preceding level.
		We combine these results to obtain the following necessary conditions.
		
		\Needspace{6\baselineskip}
		\begin{theorem}\label{the:necessary}
			Let $G$ be a connected non-bipartite $2$-threshold graph.
			Then there exist a BFS root $u$ and decompositions
			$L_i=A_i\mathbin{\dot\cup}Q_i$ for $0\le i\le d$
			for which the following conditions hold.
			
			$(1)$ For each $0\le i<d$, the family
			$\{N_{L_{i+1}}(v):v\in L_i\}$ is nested.
			
			$(2)$ Each decomposition $L_i=A_i\mathbin{\dot\cup}Q_i$ is
			admissible, and $Q_i=\emptyset$ for $0\le i\le d-2$.
			Moreover,
			$Q_{d-1}\cup N_{A_{d-1}}(Q_{d-1})\subseteq D_{d-1,s_{d-1}}$.
			
			$(3)$ For $1\le i\le d-2$, take $x_1\in D_{i,t}$ and
			$y_1\in D_{i,t'}$ with $0\le t<t'\le s_i$.
			Then $N_{L_{i-1}}(y_1)\subseteq N_{L_{i-1}}(x_1)$.
			For $i\in\{d-1,d\}$ with $i\ge1$, take $x_2\in S_{i,t}$ and
			$y_2\in S_{i,t'}$ with $0\le t<t'\le m_i$.
			Then $N_{L_{i-1}}(y_2)\subseteq N_{L_{i-1}}(x_2)$.
			
			$(4)$ For $i\in\{d-1,d\}$ with $i\ge1$,
			we have $N_{L_{i-1}}(y_3)\subseteq N_{L_{i-1}}(x_3)$
			whenever $x_3\in A_i$ and $y_3\in Q_i$.
			
			$(5)$ For $i\in\{d-1,d\}$ with $i\ge1$, take
			$x_4\in R^{(i)}_{t'}$ and $y_4\in R^{(i)}_t$, where
			$s_{d-1}\le t<t'\le m_{d-1}$ if $i=d-1$, and
			$1\le t<t'\le m_d$ if $i=d$.
			Then $N_{L_{i-1}}(y_4)\subseteq N_{L_{i-1}}(x_4)$.
		\end{theorem}
		
		\noindent\textbf{Proof.} Choose a strict representation by Lemma~$\ref{lem:strict}$.
		Take a vertex $u$ of maximum rank as the root.
		Order each odd level by increasing rank and each even level by decreasing rank.
		Corollary~$\ref{cor:next-nesting}$\,(ii) gives $(1)$.
		
		For $(2)$, Lemma~$\ref{lem:same-level}$\,(i) and~(iii) give
		admissible decompositions with $Q_i=\emptyset$ for $0\le i\le d-2$,
		and (iv) gives
		$Q_{d-1}\cup N_{A_{d-1}}(Q_{d-1})\subseteq D_{d-1,s_{d-1}}$.
		Use the classes in Definitions~$\ref{def:next-classes}$ and~$\ref{def:stable-classes}$.
		
		For $(3)$, take $x_1\in D_{i,t}$ and $y_1\in D_{i,t'}$ with
		$1\le i\le d-2$ and $0\le t<t'\le s_i$.
		Definition~$\ref{def:next-classes}$ gives
		$N_{L_{i+1}}(x_1)\subsetneqq N_{L_{i+1}}(y_1)$, so $x_1$ precedes $y_1$.
		Corollary~$\ref{cor:next-nesting}$\,(ii) then gives
		$N_{L_{i-1}}(y_1)\subseteq N_{L_{i-1}}(x_1)$.
		
		For $i\in\{d-1,d\}$ with $i\ge1$, take $x_2\in S_{i,t}$ and
		$y_2\in S_{i,t'}$ with $0\le t<t'\le m_i$.
		If $i=d-1$ and $t<s_{d-1}$, the order of the classes
		$D_{d-1,0},\ldots,D_{d-1,s_{d-1}}$ puts $x_2$ before $y_2$.
		Otherwise, suppose that $y_2$ precedes $x_2$.
		Since $A_i$ is stable, Definition~$\ref{def:stable-classes}$ and
		Lemma~$\ref{lem:same-level}$\,(ii) give
		$N_{L_i}(x_2)\subsetneqq N_{L_i}(y_2)\subseteq N_{L_i}(x_2)$,
		a contradiction.
		Thus $x_2$ precedes $y_2$ in both cases, and
		Corollary~$\ref{cor:next-nesting}$\,(ii) gives
		$N_{L_{i-1}}(y_2)\subseteq N_{L_{i-1}}(x_2)$.
		
		Conditions $(4)$ and $(5)$ follow from
		Lemma~$\ref{lem:clique-intervals}$\,(i) and~(ii), respectively. $\Box$
		
		To see why condition $(4)$ is needed, let $H$ consist of the four
		internally disjoint $a$--$b$ paths $ab$, $apb$, $aqb$, and $axyb$.
		With root $x$, we have $L_1=\{a,y\}$, $L_2=\{b,p,q\}$, and
		$N_{L_2}(y)=\{b\}\subsetneqq N_{L_2}(a)=L_2$.
		Since $L_1$ is stable, $Q_1=\emptyset$. In $L_2$, the vertices $p,q$
		are nonadjacent and each has only $b$ as a neighbor.
		Definition~$\ref{def:admissible}$ therefore forces $A_2=\{p,q\}$
		and $Q_2=\{b\}$. Conditions $(1)$ and $(2)$ hold.
		For $(3)$, both vertices of $L_1$ have preceding-level neighborhood
		$\{x\}$, and $p,q$ belong to the same class in $L_2$ since
		$N_{L_2}(p)=N_{L_2}(q)=\{b\}$. Condition $(5)$ holds since
		$Q_1=\emptyset$ and $|Q_2|=1$. Condition $(4)$ fails because
		\[
		N_{L_1}(b)=\{a,y\}\not\subseteq\{a\}=N_{L_1}(p).
		\]
		
		By symmetry, it remains to consider roots $p$ and $a$.
		With root $p$, the next-level neighborhoods of $a$ and $b$ are
		$\{q,x\}$ and $\{q,y\}$, so $(1)$ fails.
		With root $a$, the edge $pb$ lies in $L_1$ and $L_2=\{y\}$.
		Condition $(2)$ puts both ends of $pb$ in $D_{1,s_1}$ and hence requires
		both to be adjacent to $y$. Since $py\notin E(H)$, $(2)$ fails. Thus no root satisfies all five conditions, and $H$ is not
		$2$-threshold by Theorem~$\ref{the:necessary}$.
		
		\section{Construction of a representation}\label{sec:construction}
		\widowpenalty=10000
		\begin{theorem}\label{the:sufficient}
			Let $G$ be a connected graph with at least two vertices. Suppose $G$ has
			a BFS root $u$ and level decompositions satisfying conditions $(1)$--$(5)$
			of Theorem~$\ref{the:necessary}$. For any $\theta_1>3$ and
			$\theta_2>2\theta_1+2$, $G$ has a $2$-threshold representation
			$(r,(\theta_1,\theta_2))$.
		\end{theorem}
		
		Throughout this section, fix $G$, $u$, and the level decompositions as in
		Theorem~$\ref{the:sufficient}$.
		Conditions $(3)$--$(5)$ compare preceding-level neighborhoods across
		stable classes, between stable and clique vertices, and between the sets
		$R_t^{(i)}$. We use these inclusions to establish the order of the partition indices
		and to assign the adjustment values.
		
		\subsection{Partitioning the BFS levels}
		For $d\ge2$, the first three parts of $L_{d-1}$ consist of stable vertices
		with no neighbor in $Q_d\cup Q_{d-1}$, with a neighbor in $Q_d$ but none
		in $Q_{d-1}$, and with a neighbor in $Q_{d-1}$, respectively.
		The fourth part is $Q_{d-1}$.
		We then define the indices $z_k^{(i)}$ successively in the preceding levels.
		We start with the indices of classes having neighbors in $Q_d$ or $Q_{d-1}$.
		Recall $I_i=\{1,\ldots,s_i\}$ and the classes $D_{i,t}$,
		$0\le t\le s_i$, from Definition~$\ref{def:next-classes}$.
		In $A_{d-1}$, use the refined stable classes $S_{d-1,t}$,
		$0\le t\le m_{d-1}$,
		from Definition~$\ref{def:stable-classes}$\,(i), where $m_{d-1}\ge s_{d-1}$.
		Condition $(2)$ of Theorem~$\ref{the:necessary}$ gives
		$N_{Q_{d-1}}(D_{d-1,t})=\emptyset$ for $0\le t<s_{d-1}$.
		For $d\ge2$, define
		\begin{align*}
			J_1^{(d-1)}
			&=\{t\in I_{d-1}:N_{Q_d}(D_{d-1,t})\ne\emptyset\},\\
			J_2^{(d-1)}
			&=\{t\in\{s_{d-1},\ldots,m_{d-1}\}:
			N_{Q_{d-1}}(S_{d-1,t})\ne\emptyset\}.
		\end{align*}
		
		For $d\ge3$, define
		$J_3^{(d-2)}=\{t\in I_{d-2}:N_{Q_{d-1}}(D_{d-2,t})\ne\emptyset\}$.
		
		\begin{definition}[Indices defining the partition]\label{def:partition-indices}
			When $d=1$, no indices are needed. For $d\ge2$, define $z_k^{(i)}$
			for $1\le i<d$ and $k\in\{1,2,3\}$ as follows.
			\begin{enumerate}
				\item[(i)] In $L_{d-1}$, set $z_2^{(d-1)}=\min\bigl(J_2^{(d-1)}\cup\{m_{d-1}+1\}\bigr)$
				and $z_3^{(d-1)}=m_{d-1}+1$.
				If $Q_d\ne\emptyset$, then $J_1^{(d-1)}\ne\emptyset$: every vertex
				of $Q_d$ has a neighbor in $L_{d-1}$, and that neighbor belongs to
				some $D_{d-1,t}$ with $t\in I_{d-1}$. Thus set
				$z_1^{(d-1)}=\min J_1^{(d-1)}$;
				if $Q_d=\emptyset$, set $z_1^{(d-1)}=z_2^{(d-1)}$.
				
				\item[(ii)] For $d\ge3$, in $L_{d-2}$ set $z_3^{(d-2)}=\min\bigl(J_3^{(d-2)}\cup\{s_{d-2}+1\}\bigr)$.
				For $k=1,2$, define
				\[
				J_k^{(d-2)}
				=\{t\in I_{d-2}:
				N_{L_{d-1}}(D_{d-2,t})\cap S_{d-1,z_k^{(d-1)}}\ne\emptyset\},
				\]
				when $z_k^{(d-1)}\le m_{d-1}$, and set $J_k^{(d-2)}=\emptyset$
				otherwise. Put $z_k^{(d-2)}=\min\bigl(J_k^{(d-2)}\cup\{s_{d-2}+1\}\bigr)$.
				
				\item[(iii)] For $d\ge4$, define successively for $i=d-3,\ldots,1$ and
				$k\in\{1,2,3\}$,
				\[
				J_k^{(i)}
				=\{t\in I_i:N_{L_{i+1}}(D_{i,t})\cap D_{i+1,z_k^{(i+1)}}\ne\emptyset\},
				\]
				when $z_k^{(i+1)}\le s_{i+1}$, and set $J_k^{(i)}=\emptyset$ otherwise.
				Put $z_k^{(i)}=\min\bigl(J_k^{(i)}\cup\{s_i+1\}\bigr)$.
			\end{enumerate}
		\end{definition}
		
		To form the index intervals in Definition~$\ref{def:base-parts}$, we first
		verify that $z_1^{(i)}\le z_2^{(i)}\le z_3^{(i)}$.
		
		For $d\ge2$, the indices satisfy
		$1\le z_1^{(d-1)}\le z_2^{(d-1)}\le z_3^{(d-1)}=m_{d-1}+1$:
		if $Q_d\ne\emptyset$, then $z_1^{(d-1)}\le s_{d-1}\le z_2^{(d-1)}$;
		otherwise $z_1^{(d-1)}=z_2^{(d-1)}$.
		
		For $d\ge3$, first let $i=d-2$.
		If $z_2^{(d-1)}\le m_{d-1}$, conditions $(3)$ and $(4)$ of
		Theorem~$\ref{the:necessary}$ give
		\[
		N_{L_{d-2}}(Q_{d-1})
		\subseteq N_{L_{d-2}}\bigl(S_{d-1,z_2^{(d-1)}}\bigr)
		\subseteq N_{L_{d-2}}\bigl(S_{d-1,z_1^{(d-1)}}\bigr).
		\]
		The second inclusion is an equality when
		$z_1^{(d-1)}=z_2^{(d-1)}$.
		By Definition~$\ref{def:partition-indices}$,
		$J_3^{(d-2)}\subseteq J_2^{(d-2)}\subseteq J_1^{(d-2)}$, and hence
		$1\le z_1^{(d-2)}\le z_2^{(d-2)}\le z_3^{(d-2)}\le s_{d-2}+1$.
		
		If $z_2^{(d-1)}=m_{d-1}+1$, then $J_2^{(d-1)}=\emptyset$,
		so $Q_{d-1}=\emptyset$ by admissibility.
		Thus $J_2^{(d-2)}=J_3^{(d-2)}=\emptyset$, giving the same
		inclusions and inequalities.
		
		\Needspace{10\baselineskip}
		For $d\ge4$ and $i=d-3,\ldots,1$, assume that
		$z_1^{(i+1)}\le z_2^{(i+1)}\le z_3^{(i+1)}$.
		Take $1\le k<k'\le3$ with $z_{k'}^{(i+1)}\le s_{i+1}$.
		By condition $(3)$ of Theorem~$\ref{the:necessary}$,
		\[
		N_{L_i}\bigl(D_{i+1,z_{k'}^{(i+1)}}\bigr)
		\subseteq N_{L_i}\bigl(D_{i+1,z_k^{(i+1)}}\bigr),
		\]
		with equality when $z_k^{(i+1)}=z_{k'}^{(i+1)}$.
		Thus $J_{k'}^{(i)}\subseteq J_k^{(i)}$.
		
		If $z_{k'}^{(i+1)}=s_{i+1}+1$, Definition~$\ref{def:partition-indices}$
		gives $J_{k'}^{(i)}=\emptyset$, so the same inclusion holds.
		Thus $J_3^{(i)}\subseteq J_2^{(i)}\subseteq J_1^{(i)}$.
		Taking minima as in Definition~$\ref{def:partition-indices}$ gives
		$1\le z_1^{(i)}\le z_2^{(i)}\le z_3^{(i)}\le s_i+1$.
		
		\Needspace{5\baselineskip}
		\begin{definition}[Vertex parts]\label{def:base-parts}
			For $1\le i<d$, put
			\(Z_{i,1}=\{0,\ldots,z_1^{(i)}-1\}\),
			\(Z_{i,2}=\{z_1^{(i)},\ldots,z_2^{(i)}-1\}\), and
			\(Z_{i,3}=\{z_2^{(i)},\ldots,z_3^{(i)}-1\}\).
			For $1\le i\le d-2$, also put $Z_{i,4}=\{z_3^{(i)},\ldots,s_i\}$.
			An index interval is empty when its lower bound exceeds its upper bound.
			Write $B_{i,\ell}$ for part $\ell$ of level $L_i$. Set
			\[
			\begin{aligned}
				B_{i,\ell}&=\bigcup\limits_{t\in Z_{i,\ell}}D_{i,t}
				&& \text{for }1\le i\le d-2,\ 1\le\ell\le4,\\
				B_{d-1,\ell}&=\bigcup\limits_{t\in Z_{d-1,\ell}}S_{d-1,t}
				&& \text{for }d\ge2,\ 1\le\ell\le3.
			\end{aligned}
			\]
			Set $B_{d-1,4}=Q_{d-1}$ when $d\ge2$.
			In $L_d$, put $B_{d,1}=A_d$, $B_{d,2}=Q_d$, and
			$B_{d,3}=B_{d,4}=\emptyset$.
		\end{definition}
		
		For $1\le i\le d$,
		$L_i=B_{i,1}\mathbin{\dot\cup}\cdots\mathbin{\dot\cup}B_{i,4}$.
		The next lemma describes the edges between these sets.
		
		\begin{lemma}\label{lem:part-adjacency}
			For the indices and parts in Definitions~$\ref{def:partition-indices}$
			and~$\ref{def:base-parts}$, the following hold.\nopagebreak
			
			(i) If $d\ge2$ and $x_1\in A_{d-1}$, then
			$x_1\in B_{d-1,1}$ if and only if
			$N_{Q_{d-1}}(x_1)=N_{Q_d}(x_1)=\emptyset$;
			$x_1\in B_{d-1,2}$ if and only if
			$N_{Q_{d-1}}(x_1)=\emptyset$ and $N_{Q_d}(x_1)\ne\emptyset$;
			and $x_1\in B_{d-1,3}$ if and only if $N_{Q_{d-1}}(x_1)\ne\emptyset$.
			
			(ii) For $d\ge2$, $A_d\subseteq N_{L_d}(x_2)$ for every
			$x_2\in B_{d-1,2}$, and $N_{L_d}(x_3)=L_d$ for every
			$x_3\in B_{d-1,3}\cup Q_{d-1}$.
			
			(iii) For $1\le i<d$, $x_4\in B_{i,\ell_1}$, and
			$y_4\in B_{i+1,\ell_2}$, where $1\le\ell_1,\ell_2\le4$,
			$\ell_2<\ell_1$ implies $x_4y_4\in E(G)$, and
			$\ell_2>\ell_1$ implies $x_4y_4\notin E(G)$.
		\end{lemma}
		
		\noindent\textbf{Proof.} For (i), conditions $(1)$ and $(2)$ of
		Theorem~$\ref{the:necessary}$ imply that
		$\{N_{L_d}(x):x\in L_{d-1}\}$ and
		$\{N_{Q_{d-1}}(x):x\in A_{d-1}\}$ are both nested.
		For $x_1\in S_{d-1,t}$,
		Definitions~$\ref{def:t-differences}$ and~$\ref{def:r-differences}$ give
		$N_{Q_{d-1}}(x_1)=\bigcup\limits_{t_0=s_{d-1}}^tR_{t_0}^{(d-1)}$ and
		$N_{L_d}(x_1)=\bigcup\limits_{j=1}^{\min(t,s_{d-1})}T_j^{(d-1)}$.
		Definition~$\ref{def:partition-indices}$ therefore gives
		$N_{Q_{d-1}}(x_1)\ne\emptyset$ if and only if $z_2^{(d-1)}\le t$.
		When $Q_d\ne\emptyset$, it also gives
		$N_{Q_d}(x_1)\ne\emptyset$ if and only if $z_1^{(d-1)}\le t$.
		Since $z_1^{(d-1)}\le z_2^{(d-1)}$,
		$x_1\in B_{d-1,1}\Longleftrightarrow t<z_1^{(d-1)}
		\Longleftrightarrow N_{Q_{d-1}}(x_1)=N_{Q_d}(x_1)=\emptyset$.
		$x_1\in B_{d-1,2}\Longleftrightarrow z_1^{(d-1)}\le t<z_2^{(d-1)}
		\Longleftrightarrow N_{Q_{d-1}}(x_1)=\emptyset\text{ and }N_{Q_d}(x_1)\ne\emptyset$.
		$x_1\in B_{d-1,3}\Longleftrightarrow z_2^{(d-1)}\le t
		\Longleftrightarrow N_{Q_{d-1}}(x_1)\ne\emptyset$.
		
		When $Q_d=\emptyset$, $z_1^{(d-1)}=z_2^{(d-1)}$.
		The equivalences for $B_{d-1,1}$ and $B_{d-1,3}$ follow
		from the definition of $z_2^{(d-1)}$.
		The equivalence for $B_{d-1,2}$ follows from
		$B_{d-1,2}=Q_d=\emptyset$.
		
		For (ii), if $x_2\in B_{d-1,2}$, choose $y_1\in N_{Q_d}(x_2)$ by (i).
		For every $y_2\in A_d$, condition $(4)$ of Theorem~$\ref{the:necessary}$ gives
		$x_2\in N_{L_{d-1}}(y_1)\subseteq N_{L_{d-1}}(y_2)$.
		Thus $A_d\subseteq N_{L_d}(x_2)$.
		By (i) and condition $(2)$ of Theorem~$\ref{the:necessary}$,
		$B_{d-1,3}\cup Q_{d-1}
		=N_{A_{d-1}}(Q_{d-1})\cup Q_{d-1}\subseteq D_{d-1,s_{d-1}}$.
		Definition~$\ref{def:next-classes}$ and Corollary~$\ref{cor:next-nesting}$\,(i) give
		$N_{L_d}(x_3)=L_d$ for every $x_3\in B_{d-1,3}\cup Q_{d-1}$.
		
		For (iii), first let $i=d-1\ge1$.
		Since $B_{d,3}=B_{d,4}=\emptyset$, we have $\ell_2\in\{1,2\}$.
		If $\ell_2<\ell_1$ and $\ell_1=2$, then $y_4\in A_d$
		and (ii) gives $x_4y_4\in E(G)$.
		If $\ell_2<\ell_1$ and $\ell_1\in\{3,4\}$, then
		$N_{L_d}(x_4)=L_d$ by (ii), so $x_4y_4\in E(G)$.
		If $\ell_2>\ell_1$, then $\ell_1=1$, $\ell_2=2$;
		(i) gives $x_4y_4\notin E(G)$.
		
		Let $1\le i\le d-2$. Conditions $(3)$ and $(4)$ of Theorem~$\ref{the:necessary}$ give
		\[
		\begin{aligned}
			N_{L_i}\!\left(\bigcup\limits_{t=t_0}^{s_{i+1}}D_{i+1,t}\right)
			&=N_{L_i}(D_{i+1,t_0})
			&&\text{for }i\le d-3,\ 1\le t_0\le s_{i+1},\\
			N_{L_{d-2}}\!\left(Q_{d-1}\cup\bigcup\limits_{t=t_0}^{m_{d-1}}S_{d-1,t}\right)
			&=N_{L_{d-2}}(S_{d-1,t_0})
			&&\text{for }1\le t_0\le m_{d-1}.
		\end{aligned}
		\]
		Consequently, Definition~$\ref{def:partition-indices}$ and condition $(1)$ of Theorem~$\ref{the:necessary}$
		give, for $k\in\{1,2,3\}$,
		\[
		N_{L_i}\!\left(\bigcup\limits_{h=k+1}^4B_{i+1,h}\right)
		=\bigcup\limits_{t\in J_k^{(i)}}D_{i,t}
		=\bigcup\limits_{t=z_k^{(i)}}^{s_i}D_{i,t}.
		\]
		When $z_k^{(i)}=s_i+1$,
		$\bigcup\limits_{t=z_k^{(i)}}^{s_i}D_{i,t}=\emptyset$.
		For $i=d-2$ and $k=3$,
		$\bigcup\limits_{h=k+1}^4B_{i+1,h}=B_{d-1,4}=Q_{d-1}$.
		
		Write $x_4\in D_{i,t}$.
		If $\ell_2<\ell_1$, then $t\ge z_{\ell_1-1}^{(i)}$,
		so taking $k=\ell_1-1$ in the identity above gives a vertex
		$v\in\bigcup\limits_{h=\ell_1}^4B_{i+1,h}$ with $x_4v\in E(G)$.
		When $y_4,v$ lie in stable classes, their index ranges and condition $(3)$
		of Theorem~$\ref{the:necessary}$ give
		$x_4\in N_{L_i}(v)\subseteq N_{L_i}(y_4)$, so $x_4y_4\in E(G)$.
		If $i=d-2$, $y_4\in A_{d-1}$, and $v\in Q_{d-1}$,
		condition $(4)$ of Theorem~$\ref{the:necessary}$ gives the same inclusion
		and hence $x_4y_4\in E(G)$.
		
		If $\ell_2>\ell_1$, then $t<z_{\ell_1}^{(i)}$ and
		$y_4\in\bigcup\limits_{h=\ell_1+1}^4B_{i+1,h}$.
		The same equality gives $x_4y_4\notin E(G)$.
		When $d=1$, (iii) is vacuous. $\Box$
		
		\subsection{The base functions}
		For the pairs in Lemma~$\ref{lem:part-adjacency}$ with different indices $\ell$,
		the base values determine whether the rank sum lies between the thresholds.
		When the two base values sum to a threshold, the adjustments determine
		on which side of that threshold the rank sum lies.
		\begin{definition}[Base values]\label{def:base-values}
			Choose $\theta_1>3$ and $\theta_2>2\theta_1+2$.
			For $0\le i\le d$, define
			\begin{equation}\label{eq:base-first}
				f_1(i)=
				\begin{cases}
					\theta_1-\frac{d-i+1}{2}(\theta_2-\theta_1),
					& \text{if $d-i$ is odd},\\[6pt]
					\frac{d-i+2}{2}(\theta_2-\theta_1),
					& \text{if $d-i$ is even},
				\end{cases}
			\end{equation}
			and, for $\ell=2,3,4$, put
			\begin{equation}\label{100}
				f_\ell(i)
				=
				f_1(i)
				+(-1)^{d-i+1}
				\left(
				\frac{(\ell-2)\theta_1}{3}+\theta_2-2\theta_1
				\right).
			\end{equation}
		\end{definition}
		
		\begin{lemma}\label{lem:base-sums}
			The base values in Definition~$\ref{def:base-values}$ satisfy
			the following identities and inequalities.
			
			(i) For $0\le i\le d$,
			$f_1(i)<f_2(i)<f_3(i)<f_4(i)$ when $d-i$ is odd, and
			$f_4(i)<f_3(i)<f_2(i)<f_1(i)$ when $d-i$ is even.
			
			(ii) For $0\le i\le d$ and $1\le p<q\le4$,
			$1<|f_q(i)-f_p(i)|<\theta_2-\theta_1-1$.
			
			(iii) For $0\le i<d$ and $1\le\ell\le4$,
			$f_\ell(i)+f_\ell(i+1)=\theta_1$ when $d-i$ is odd and
			$f_\ell(i)+f_\ell(i+1)=\theta_2$ when $d-i$ is even.
			
			(iv) At the last two levels,
			$f_3(d-1)+f_4(d-1)=\theta_1$ and $f_1(d)+f_2(d)=\theta_2$.
		\end{lemma}
		
		\noindent\textbf{Proof.} For (i) and (ii), we distinguish the parity of $d-i$.
		For $d-i=2k+1$ with $k\ge0$,
		equations~\eqref{eq:base-first} and~\eqref{100} give
		$f_1(i)=\theta_1-(k+1)(\theta_2-\theta_1)$ and
		$f_\ell(i)=-k(\theta_2-\theta_1)+\frac{(\ell-2)\theta_1}{3}$
		for $\ell\in\{2,3,4\}$. Hence
		$f_1(i)<f_2(i)<f_3(i)<f_4(i)$.
		Here $f_2(i)-f_1(i)=\theta_2-2\theta_1>2$ and
		$f_3(i)-f_2(i)=f_4(i)-f_3(i)=\frac{\theta_1}{3}>1$. Therefore
		$1<|f_q(i)-f_p(i)|\le|f_4(i)-f_1(i)|
		=\theta_2-\frac{4\theta_1}{3}<\theta_2-\theta_1-1$,
		where the last inequality uses $\frac{\theta_1}{3}>1$.
		
		For $d-i=2k$ with $k\ge0$, equations~\eqref{eq:base-first} and~\eqref{100} give
		$f_1(i)=(k+1)(\theta_2-\theta_1)$ and
		$f_\ell(i)=k(\theta_2-\theta_1)+\frac{(5-\ell)\theta_1}{3}$
		for $\ell\in\{2,3,4\}$. Thus
		$f_4(i)<f_3(i)<f_2(i)<f_1(i)$.
		Here $f_1(i)-f_2(i)=\theta_2-2\theta_1>2$ and
		$f_2(i)-f_3(i)=f_3(i)-f_4(i)=\frac{\theta_1}{3}>1$. Therefore
		$1<|f_q(i)-f_p(i)|\le\theta_2-\frac{4\theta_1}{3}
		<\theta_2-\theta_1-1$ for $p<q$.
		
		For (iii), let $0\le i<d$. For $\ell\in\{2,3,4\}$,
		equation~\eqref{100} gives
		$f_\ell(i)-f_1(i)
		=-\bigl(f_\ell(i+1)-f_1(i+1)\bigr)$.
		Hence $f_\ell(i)+f_\ell(i+1)=f_1(i)+f_1(i+1)$ for $1\le\ell\le4$.
		If $d-i=2k+1$,
		equation~\eqref{eq:base-first} gives
		$f_1(i)+f_1(i+1)
		=\theta_1-(k+1)(\theta_2-\theta_1)+(k+1)(\theta_2-\theta_1)=\theta_1$.
		Thus $f_\ell(i)+f_\ell(i+1)=\theta_1$ for every $1\le\ell\le4$.
		For $d-i=2k$ with $k\ge1$,
		$f_1(i)+f_1(i+1)
		=(k+1)(\theta_2-\theta_1)+\theta_1-k(\theta_2-\theta_1)=\theta_2$.
		Thus $f_\ell(i)+f_\ell(i+1)=\theta_2$
		for every $1\le\ell\le4$.
		
		For (iv), in $L_{d-1}$, equations~\eqref{eq:base-first} and~\eqref{100}
		give $f_3(d-1)+f_4(d-1)=\theta_1$,
		since $f_3(d-1)=\frac{\theta_1}{3}$ and $f_4(d-1)=\frac{2\theta_1}{3}$.
		In the last level, $f_1(d)+f_2(d)=\theta_2$,
		since $f_1(d)=\theta_2-\theta_1$ and $f_2(d)=\theta_1$. $\Box$
		
		\subsection{The adjustment and rank functions}
		\begin{definition}[Adjustment function]\label{def:adjustment-values}
			Let $n=|V(G)|$. For $1\le i\le d$, $1\le\ell\le4$, and
			$v\in B_{i,\ell}$, let $t(v)$ be the unique class index $t$ such that
			$v\in D_{i,t}$ when $i\le d-2$, $v\in S_{i,t}$ when
			$v\in A_i$ and $i\in\{d-1,d\}$, or $v\in R_t^{(i)}$ when $v\in Q_i$.
			Let $j(v)$ be the unique class index $j\in\{1,\ldots,s_{i-1}\}$ such that
			$v\in T_j^{(i-1)}$.
			Here $0\le t(v)\le s_i$ when $i\le d-2$, and
			$0\le t(v)\le m_i$ when $i\in\{d-1,d\}$ and $v\in A_i$.
			For $v\in Q_{d-1}$ and $v\in Q_d$, the respective ranges are
			$s_{d-1}\le t(v)\le m_{d-1}$ and $1\le t(v)\le m_d$.
			For each nonempty $B_{i,\ell}$, define the coefficient function of the class indices $t$ and $j$ by
			\begin{equation}\label{eq:adjustment-parameter}
				\eta_{i,\ell}(t,j)=
				\begin{cases}
					1-\dfrac{nt+j}{n^2},&\ell\in\{1,3\},\ B_{i,\ell}\subseteq A_i,\\[5pt]
					\dfrac{nt+j}{n^2},&\ell\in\{2,4\},\ B_{i,\ell}\subseteq A_i,\\[5pt]
					1-\dfrac{nt-j}{n^2},&B_{i,\ell}\subseteq Q_i.
				\end{cases}
			\end{equation}
			For $v\in B_{i,\ell}$, we use the coefficient
			$\eta_{i,\ell}(t(v),j(v))$.
			
			For integer indices $t<t'$ and $1\le j,j'<n$,
			$nt+j<nt'+j'$ and $nt-j<nt'-j'$.
			
			\Needspace{5\baselineskip}
			We define $\phi:V(G)\setminus\{u\}\to(0,\tfrac12)$ successively,
			from left to right in each of the following two sequences:
			\[
			\begin{aligned}
				&B_{1,1},\ldots,B_{d,1},B_{d,2},\ldots,B_{1,2},\\
				&B_{1,3},\ldots,B_{d-1,3},B_{d-1,4},\ldots,B_{1,4}.
			\end{aligned}
			\]
			For $d=1$, omit the second sequence.
			These orders make each pair of nonempty parts compared in
			Lemma~$\ref{lem:adjustment-adjacency}$ consecutive.
			In each displayed sequence, let $B^-_{i,\ell}$ be the set immediately
			preceding $B_{i,\ell}$.
			Set $B^-_{1,1}=\emptyset$ and, when $d\ge2$, $B^-_{1,3}=\emptyset$;
			thus $B^-_{d,2}=B_{d,1}$.
			For $v,w\in B_{i,\ell}$, Claim~\ref{clm:current-neighborhoods} shows that
			$\eta_{i,\ell}(t(v),j(v))<\eta_{i,\ell}(t(w),j(w))$
			implies $N_{B^-_{i,\ell}}(v)\subseteq N_{B^-_{i,\ell}}(w)$.
			For $v\in B_{i,\ell}$, only values of $\phi$ already assigned on
			$B^-_{i,\ell}$ are used in the following definition.
			If $B^-_{i,\ell}=\emptyset$, set
			$\phi_-(v)=0$, $\phi_+(v)=\tfrac12$, and
			$\phi(v)=\tfrac12\eta_{i,\ell}\bigl(t(v),j(v)\bigr)$.
			Otherwise, the values on $B^-_{i,\ell}$ have already been defined; put
			\[
			\begin{aligned}
				\phi_-(v)&=\max\bigl(\{\phi(x):x\in B^-_{i,\ell}\cap N_G(v)\}\cup\{0\}\bigr),\\
				\phi_+(v)&=\min\bigl(\{\phi(x):x\in B^-_{i,\ell}\setminus N_G(v)\}\cup\{\tfrac12\}\bigr),\\
				\phi(v)&=\phi_-(v)+\eta_{i,\ell}\bigl(t(v),j(v)\bigr)
				\bigl(\phi_+(v)-\phi_-(v)\bigr).
			\end{aligned}
			\]
		\end{definition}
		
		For the root $u\in L_0$, choose the base value $f_4(0)$ and set
		$r(u)=f_4(0)+(-1)^{d+1}\frac12$.
		For $1\le i\le d$ and $v\in L_i$, define
		\[
		r(v)=
		\begin{cases}
			f_1(i)+(-1)^{d-i}\phi(v),&v\in B_{i,1},\\[4pt]
			f_2(i)+(-1)^{d-i+1}\phi(v),&v\in B_{i,2},\\[4pt]
			f_3(i)+(-1)^{d-i}\phi(v),&v\in B_{i,3},\\[4pt]
			f_4(i)+(-1)^{d-i+1}\phi(v),&v\in B_{i,4}.
		\end{cases}
		\]
		\Needspace{4\baselineskip}
		\begin{lemma}\label{lem:adjustments}
			The function $\phi$ in Definition~$\ref{def:adjustment-values}$ is well defined.
			For all $v,w\in V(G)\setminus\{u\}$, the following hold.
			
			(i) $\phi(v)\in\mathbb Q$.
			
			(ii) $0<\phi(v)<\tfrac12$.
			
			(iii) If $v,w\in B_{i,\ell}$ satisfy $t(v)=t(w)$ and $j(v)=j(w)$,
			then $\phi(v)=\phi(w)$.
			
			(iv) $\phi_-(v)<\phi(v)<\phi_+(v)$.
			
			(v) If $v,w\in B_{i,\ell}$ and
			$\eta_{i,\ell}(t(v),j(v))<\eta_{i,\ell}(t(w),j(w))$,
			then~${\phi(v)<\phi(w)}$.
		\end{lemma}
		
		\noindent\textbf{Proof.} For (i), since $n,t(v),j(v)$ are integers,
		equation~\eqref{eq:adjustment-parameter} gives
		$\eta_{i,\ell}(t(v),j(v))\in\mathbb Q$.
		We first consider $B_{1,1},\ldots,B_{d,1},B_{d,2},\ldots,B_{1,2}$.
		For $v\in B_{1,1}$, the identity $B^-_{1,1}=\emptyset$ gives
		$\phi(v)=\tfrac12\eta_{1,1}(t(v),j(v))\in\mathbb Q$.
		
		For $2\le i\le d$, assume that $\phi(x)\in\mathbb Q$
		for every $x\in B_{i-1,1}$.
		Since $B^-_{i,1}=B_{i-1,1}$, the definitions give
		$\phi_-(v),\phi_+(v)\in\{\phi(x):x\in B_{i-1,1}\}\cup\{0,\tfrac12\}\subseteq\mathbb Q$
		for every $v\in B_{i,1}$.
		Thus $\phi(v)=\phi_-(v)+\eta_{i,1}(t(v),j(v))\bigl(\phi_+(v)-\phi_-(v)\bigr)\in\mathbb Q$.
		This proves (i) on $B_{1,1},\ldots,B_{d,1}$.
		
		Next, $B^-_{d,2}=B_{d,1}$ gives (i) on $B_{d,2}$ by the same argument.
		When $d\ge2$, the identities $B^-_{i,2}=B_{i+1,2}$ for $1\le i<d$
		give (i) successively on $B_{d-1,2},\ldots,B_{1,2}$.
		The second sequence is treated similarly.
		
		For (ii), we first bound $t(v)$ and $j(v)$ to show
		$0<\eta_{i,\ell}(t(v),j(v))<1$. Fix $v\in B_{i,\ell}$.
		If $i\le d-2$, then $v\in D_{i,t(v)}$.
		The sets $D_{i,1},\ldots,D_{i,s_i}$ are nonempty and pairwise disjoint
		subsets of $L_i$, so $0\le t(v)\le s_i\le|L_i|<n$.
		
		If $i\in\{d-1,d\}$, then $v\in S_{i,t(v)}$ or $v\in R_{t(v)}^{(i)}$,
		and both definitions give $t(v)\le m_i$.
		The sets $S_{i,t}$ with $1\le t\le m_i$ are nonempty and pairwise disjoint
		subsets of $A_i$, so $0\le t(v)\le m_i\le|A_i|<n$.
		
		Also, $v\in T_{j(v)}^{(i-1)}$, where the classes $T_j^{(i-1)}$
		are indexed by $1\le j\le s_{i-1}$.
		The sets $D_{i-1,1},\ldots,D_{i-1,s_{i-1}}$ are nonempty and pairwise disjoint
		subsets of $L_{i-1}$, so $1\le j(v)\le s_{i-1}\le|L_{i-1}|<n$.
		
		Since $0\le t(v)\le n-1$ and $1\le j(v)\le n-1$,
		$1\le nt(v)+j(v)\le n(n-1)+(n-1)=n^2-1<n^2$.
		For $v\in Q_{d-1}\cup Q_d$, we have $t(v)\ge1$, so
		$0<nt(v)-j(v)<n^2$.
		Equation~\eqref{eq:adjustment-parameter} gives
		$0<\eta_{i,\ell}(t(v),j(v))<1$ for $v\in B_{i,\ell}$.
		
		\Needspace{4\baselineskip}
		If $B^-_{i,\ell}=\emptyset$, then
		$0<\phi(v)=\tfrac12\eta_{i,\ell}(t(v),j(v))<\tfrac12$.
		If $B^-_{i,\ell}\ne\emptyset$ and $0<\phi(x)<\tfrac12$
		for every $x\in B^-_{i,\ell}$, the definitions give
		$0\le\phi_-(v)<\tfrac12$ and $0<\phi_+(v)\le\tfrac12$.
		Hence $0<\phi(v)=\bigl(1-\eta_{i,\ell}(t(v),j(v))\bigr)\phi_-(v)
		+\eta_{i,\ell}(t(v),j(v))\phi_+(v)<\tfrac12$.
		
		For (iii), take $v,w\in B_{i,\ell}$ with $t(v)=t(w)=t$ and $j(v)=j(w)=j$.
		The class definitions and conditions $(1)$ and $(2)$ of Theorem~$\ref{the:necessary}$ give
		$N_{B^-_{i,\ell}}(v)=N_{B^-_{i,\ell}}(w)$, so
		$\phi_-(v)=\phi_-(w)$ and $\phi_+(v)=\phi_+(w)$.
		Thus $\phi(v)=\phi_-(w)+\eta_{i,\ell}(t,j)\bigl(\phi_+(w)-\phi_-(w)\bigr)=\phi(w)$.
		
		For (iv) and (v), we first consider
		$B_{1,1},\ldots,B_{d,1},B_{d,2},\ldots,B_{1,2}$.
		Since $B^-_{1,1}=\emptyset$, for $y_4\in B_{1,1}$ the definition gives
		$\phi_-(y_4)=0<\phi(y_4)=\tfrac12\eta_{1,1}(t(y_4),j(y_4))<\tfrac12=\phi_+(y_4)$.
		For $y_5\in B_{1,1}$ with $\eta_{1,1}(t(y_4),j(y_4))<\eta_{1,1}(t(y_5),j(y_5))$,
		$\phi(y_5)-\phi(y_4)=\tfrac12\bigl(\eta_{1,1}(t(y_5),j(y_5))-\eta_{1,1}(t(y_4),j(y_4))\bigr)>0$.
		If $B^-_{i,\ell}=\emptyset$ with $\ell\in\{1,2\}$,
		the same formulas prove (iv) and (v) on $B_{i,\ell}$.
		
		Suppose (iv) and (v) hold on a nonempty set $B=B_{i',\ell'}$ in the sequence.
		We prove them on the next set $B_{i,\ell}$ when it is nonempty.
		To show that $\phi_-(y_6)<\phi_+(y_6)$ for $y_6\in B_{i,\ell}$,
		we use the next claim to prove $\phi(x_4)<\phi(x_5)$ whenever
		$x_4\in N_B(y_6)$ and $x_5\in B\setminus N_B(y_6)$.
		\Needspace{6\baselineskip}
		\begin{claim}\label{clm:predecessor-neighborhoods}
			For $x_2,x_3\in B$ with $\eta_{i',\ell'}(t(x_2),j(x_2))\le\eta_{i',\ell'}(t(x_3),j(x_3))$,
			$N_{B_{i,\ell}}(x_3)\subseteq N_{B_{i,\ell}}(x_2)$.
		\end{claim}
		
		\noindent\textit{Proof.}
		If $B_{i,\ell}\subseteq A_i$ and $\ell\in\{1,3\}$, then
		$B=B_{i-1,\ell}$. By equation~\eqref{eq:adjustment-parameter} and $|j(x_2)-j(x_3)|<n$,
		$t(x_2)-t(x_3)\ge \tfrac{j(x_3)-j(x_2)}{n}>-1$.
		Since $t(x_2)-t(x_3)$ is an integer, $t(x_2)\ge t(x_3)$.
		Let $z\in B_{i,\ell}$. By Definition~$\ref{def:adjustment-values}$,
		$z\in T^{(i-1)}_{j(z)}$. If $i<d$, we also have $x_2\in D_{i-1,t(x_2)}$,
		so Corollary~$\ref{cor:t-adjacency}$\,(ii) gives
		$x_2z\in E(G)\Longleftrightarrow j(z)\le t(x_2)$.
		The same equivalence holds for $x_3$, so $t(x_2)\ge t(x_3)$ gives
		\[
		N_{B_{i,\ell}}(x_3)
		=\{z\in B_{i,\ell}:j(z)\le t(x_3)\}
		\subseteq\{z\in B_{i,\ell}:j(z)\le t(x_2)\}
		=N_{B_{i,\ell}}(x_2).
		\]
		For $i=d$, $x_2\in S_{d-1,t(x_2)}\subseteq
		D_{d-1,\min\{t(x_2),s_{d-1}\}}$.
		By Corollary~$\ref{cor:t-adjacency}$\,(ii) and $j(z)\le s_{d-1}$,
		$x_2z\in E(G)\Longleftrightarrow j(z)\le\min\{t(x_2),s_{d-1}\}
		\Longleftrightarrow j(z)\le t(x_2)$.
		The same equivalence holds for $x_3$, so the inclusion follows
		as in the case $i<d$.
		
		When $B_{i,\ell}\subseteq A_i$ and $\ell\in\{2,4\}$, the preceding set is
		$B=B_{i+1,\ell}$. If $B\subseteq A_{i+1}$, equation~\eqref{eq:adjustment-parameter} gives
		$n\bigl(t(x_2)-t(x_3)\bigr)\le j(x_3)-j(x_2)$.
		Since $|j(x_2)-j(x_3)|<n$,
		$t(x_2)\le t(x_3)$.
		When $t(x_2)<t(x_3)$, condition $(3)$ of Theorem~$\ref{the:necessary}$ gives
		$N_{L_i}(x_3)\subseteq N_{L_i}(x_2)$, and hence
		$N_{B_{i,\ell}}(x_3)\subseteq N_{B_{i,\ell}}(x_2)$.
		When $t(x_2)=t(x_3)$, equation~\eqref{eq:adjustment-parameter} gives
		$j(x_2)\le j(x_3)$. By Corollary~$\ref{cor:t-adjacency}$\,(ii),
		\[
		N_{B_{i,\ell}}(x_3)=B_{i,\ell}\cap\bigcup\limits_{h=j(x_3)}^{s_i}D_{i,h}
		\subseteq B_{i,\ell}\cap\bigcup\limits_{h=j(x_2)}^{s_i}D_{i,h}=N_{B_{i,\ell}}(x_2).
		\]
		
		If $B\not\subseteq A_{i+1}$, Definition~$\ref{def:base-parts}$ gives
		$B=Q_{i+1}$ with $i+1\in\{d-1,d\}$.
		Equation~\eqref{eq:adjustment-parameter} gives
		$n\bigl(t(x_2)-t(x_3)\bigr)\ge j(x_2)-j(x_3)$.
		Since $|j(x_2)-j(x_3)|<n$,
		$t(x_2)\ge t(x_3)$.
		When $t(x_2)>t(x_3)$, condition $(5)$ of Theorem~$\ref{the:necessary}$ gives
		$N_{L_i}(x_3)\subseteq N_{L_i}(x_2)$.
		When $t(x_2)=t(x_3)$, equation~\eqref{eq:adjustment-parameter} gives
		$j(x_2)\le j(x_3)$, so $N_{L_i}(x_3)\subseteq N_{L_i}(x_2)$ by
		Corollary~$\ref{cor:t-adjacency}$\,(iii).
		Since $B_{i,\ell}\subseteq L_i$, both cases give
		\[
		N_{B_{i,\ell}}(x_3)
		= N_{L_i}(x_3)\cap B_{i,\ell}
		\subseteq N_{L_i}(x_2)\cap B_{i,\ell}
		= N_{B_{i,\ell}}(x_2).
		\]
		
		If $B_{i,\ell}=Q_i$ and $i=d-1\ge1$, then $B=B_{d-1,3}$.
		Equation~\eqref{eq:adjustment-parameter} gives
		$n\bigl(t(x_2)-t(x_3)\bigr)\ge j(x_3)-j(x_2)$.
		Since $|j(x_2)-j(x_3)|<n$,
		$t(x_2)\ge t(x_3)$.
		By Corollary~$\ref{cor:r-adjacency}$\,(i)--(ii),
		\[
		N_{Q_{d-1}}(x_3)=\bigcup\limits_{h=s_{d-1}}^{t(x_3)}R_h^{(d-1)}
		\subseteq\bigcup\limits_{h=s_{d-1}}^{t(x_2)}R_h^{(d-1)}=N_{Q_{d-1}}(x_2).
		\]
		
		If $B_{i,\ell}=Q_i$ and $i=d$, then $B=B_{d,1}$.
		Analogously to the case $i=d-1$, equation~\eqref{eq:adjustment-parameter}
		gives $t(x_2)\ge t(x_3)$.
		By Corollary~$\ref{cor:r-adjacency}$\,(i)--(ii),
		\[
		N_{Q_d}(x_3)=\bigcup\limits_{h=1}^{t(x_3)}R_h^{(d)}
		\subseteq\bigcup\limits_{h=1}^{t(x_2)}R_h^{(d)}=N_{Q_d}(x_2).
		\]
		This proves the claim. $\Box$
		
		For (iv), fix $y_6\in B_{i,\ell}$.
		Suppose $\emptyset\subsetneqq N_B(y_6)\subsetneqq B$, and take arbitrary
		$x_4\in N_B(y_6)$ and $x_5\in B\setminus N_B(y_6)$.
		If $\eta_{i',\ell'}(t(x_5),j(x_5))\le\eta_{i',\ell'}(t(x_4),j(x_4))$, Claim~\ref{clm:predecessor-neighborhoods} gives
		$y_6\in N_{B_{i,\ell}}(x_4)\subseteq N_{B_{i,\ell}}(x_5)$,
		contrary to $x_5\notin N_B(y_6)$. Consequently,
		$\eta_{i',\ell'}(t(x_4),j(x_4))<\eta_{i',\ell'}(t(x_5),j(x_5))$,
		and $\phi(x_4)<\phi(x_5)$ by (v) on $B$.
		By (ii) and the definitions of $\phi_-$ and $\phi_+$,
		$0<\phi(x_4)\le\phi_-(y_6)$ and
		$\phi_+(y_6)\le\phi(x_5)<\tfrac12$.
		Since $x_4$ and $x_5$ are arbitrary, we obtain
		\[
		\phi_-(y_6)=\max\limits_{x\in N_B(y_6)}\phi(x)
		<\min\limits_{x\in B\setminus N_B(y_6)}\phi(x)=\phi_+(y_6).
		\]
		
		If $N_B(y_6)=\emptyset$, (ii) and Definition~$\ref{def:adjustment-values}$ give
		$\phi_-(y_6)=0<\min\limits_{x\in B}\phi(x)=\phi_+(y_6)$.
		If $N_B(y_6)=B$, then
		$\phi_-(y_6)=\max\limits_{x\in B}\phi(x)<\tfrac12=\phi_+(y_6)$.
		
		Thus $\phi_-(y_6)<\phi_+(y_6)$. Since
		$0<\eta_{i,\ell}(t(y_6),j(y_6))<1$ by the proof of (ii),
		\[
		\phi_-(y_6)<\phi(y_6)
		=\phi_-(y_6)+\eta_{i,\ell}\bigl(t(y_6),j(y_6)\bigr)
		\bigl(\phi_+(y_6)-\phi_-(y_6)\bigr)<\phi_+(y_6).
		\]
		
		To prove $\phi(y_6)<\phi(y_7)$ in (v), we first show that
		$N_B(y_6)\subseteq N_B(y_7)$.
		\begin{claim}\label{clm:current-neighborhoods}
			For $y_6,y_7\in B_{i,\ell}$ with
			$\eta_{i,\ell}(t(y_6),j(y_6))<\eta_{i,\ell}(t(y_7),j(y_7))$,
			$N_B(y_6)\subseteq N_B(y_7)$.
		\end{claim}
		
		\noindent\textit{Proof.}
		If $B_{i,\ell}\subseteq A_i$ and $\ell\in\{1,3\}$, then
		$B=B_{i-1,\ell}$. As in Claim~\ref{clm:predecessor-neighborhoods}, equation~\eqref{eq:adjustment-parameter} gives
		$t(y_6)\ge t(y_7)$.
		When $t(y_6)>t(y_7)$, condition $(3)$ of Theorem~$\ref{the:necessary}$ gives
		$N_B(y_6)\subseteq N_B(y_7)$.
		When $t(y_6)=t(y_7)$, equation~\eqref{eq:adjustment-parameter} gives
		$j(y_6)>j(y_7)$. By Corollary~$\ref{cor:t-adjacency}$\,(ii) for $i<d$
		and~(iv) for $i=d$,
		\[
		N_B(y_6)=\{z\in B:j(y_6)\le t(z)\}
		\subseteq\{z\in B:j(y_7)\le t(z)\}=N_B(y_7).
		\]
		
		If $B_{i,\ell}\subseteq A_i$ and $\ell\in\{2,4\}$, then
		$B=B_{i+1,\ell}$ and equation~\eqref{eq:adjustment-parameter} gives
		$t(y_6)\le t(y_7)$.
		By Corollary~$\ref{cor:t-adjacency}$\,(ii) for $i\le d-2$
		and~(iv) for $i=d-1$,
		\[
		N_B(y_6)=\{z\in B:j(z)\le t(y_6)\}
		\subseteq\{z\in B:j(z)\le t(y_7)\}=N_B(y_7).
		\]
		
		If $B_{i,\ell}=Q_i$, equation~\eqref{eq:adjustment-parameter} gives
		$t(y_6)\ge t(y_7)$.
		Since $B\subseteq A_i$, Corollary~$\ref{cor:r-adjacency}$\,(ii) gives
		\[
		N_B(y_6)=\{z\in B:t(y_6)\le t(z)\}
		\subseteq\{z\in B:t(y_7)\le t(z)\}=N_B(y_7).
		\tag*{$\Box$}
		\]
		
		Now take $y_7\in B_{i,\ell}$ with
		$\eta_{i,\ell}(t(y_6),j(y_6))<\eta_{i,\ell}(t(y_7),j(y_7))$.
		Claim~\ref{clm:current-neighborhoods} gives $N_B(y_6)\subseteq N_B(y_7)$.
		If $N_B(y_6)=N_B(y_7)$, the definitions give the same $\phi_-$ and $\phi_+$, so
		\[
		\phi(y_7)-\phi(y_6)
		=(\eta_{i,\ell}(t(y_7),j(y_7))-\eta_{i,\ell}(t(y_6),j(y_6)))\bigl(\phi_+(y_6)-\phi_-(y_6)\bigr)>0.
		\]
		If $N_B(y_6)\subsetneqq N_B(y_7)$, choose
		$x_6\in N_B(y_7)\setminus N_B(y_6)$.
		The definitions of $\phi_+$ and $\phi_-$ give
		\[
		\phi(y_6)<\phi_+(y_6)\le\phi(x_6)\le\phi_-(y_7)<\phi(y_7).
		\]
		Thus (iv) and (v) hold on $B_{i,\ell}$.
		With $\ell=1$ and $B=B_{i-1,1}$ for $2\le i\le d$,
		this proves (iv) and (v) on $B_{1,1},\ldots,B_{d,1}$.
		Next, $B^-_{d,2}=B_{d,1}$ gives (iv) and (v) on $B_{d,2}$ by the same argument.
		When $d\ge2$, the identities $B^-_{i,2}=B_{i+1,2}$ for $1\le i<d$
		give (iv) and (v) successively on $B_{d-1,2},\ldots,B_{1,2}$.
		The second sequence is treated similarly. $\Box$
		
		\Needspace{4\baselineskip}
		\begin{lemma}\label{lem:adjustment-adjacency}
			For the adjustment function $\phi$ in Definition~$\ref{def:adjustment-values}$,
			the following hold.
			
			(i) For $1\le i<d$, $x\in B_{i,\ell}$, and $y\in B_{i+1,\ell}$,
			if $\ell\in\{1,3\}$, then $xy\in E(G)$ if and only if $\phi(x)<\phi(y)$;
			if $\ell\in\{2,4\}$, then $xy\in E(G)$ if and only if $\phi(x)>\phi(y)$.
			For a nonedge, $\phi(x)>\phi(y)$ when $\ell\in\{1,3\}$,
			and $\phi(x)<\phi(y)$ when $\ell\in\{2,4\}$.
			
			(ii) For $x\in A_d$ and $y\in Q_d$, and also for
			$x\in B_{d-1,3}$ and $y\in Q_{d-1}$ when $d\ge2$,
			$xy\in E(G)$ if and only if $\phi(x)<\phi(y)$.
			For a nonedge, $\phi(x)>\phi(y)$.
		\end{lemma}
		
		\noindent\textbf{Proof.} For (i), take $x\in B_{i,\ell}$ and
		$y\in B_{i+1,\ell}$ with $1\le i<d$.
		Suppose $\ell\in\{1,3\}$. If $\ell=3$, then $i\le d-2$ because
		$B_{d,3}=\emptyset$.
		For fixed $\ell$, Definition~$\ref{def:adjustment-values}$ lists the sets
		$B_{i,\ell}$ in increasing order of $i$, so $B^-_{i+1,\ell}=B_{i,\ell}$.
		If $xy\in E(G)$, then $x\in B_{i,\ell}\cap N_G(y)$.
		By the definition of $\phi_-$ and Lemma~$\ref{lem:adjustments}$\,(iv),
		\[
		\phi(x)\le
		\max\bigl(\{\phi(z):z\in B_{i,\ell}\cap N_G(y)\}\cup\{0\}\bigr)
		=\phi_-(y)<\phi(y).
		\]
		If $xy\notin E(G)$, then $x\in B_{i,\ell}\setminus N_G(y)$.
		By Lemma~$\ref{lem:adjustments}$\,(iv) and the definition of $\phi_+$,
		\[
		\phi(y)<\phi_+(y)
		=\min\bigl(\{\phi(z):z\in B_{i,\ell}\setminus N_G(y)\}\cup\{\tfrac12\}\bigr)
		\le\phi(x).
		\]
		Thus $xy\in E(G)\Longleftrightarrow\phi(x)<\phi(y)$.
		
		For $\ell\in\{2,4\}$, $i\le d-2$ when $\ell=4$ because
		$B_{d,4}=\emptyset$.
		For fixed $\ell$, Definition~$\ref{def:adjustment-values}$ lists the sets
		$B_{i,\ell}$ in decreasing order of $i$, so $B^-_{i,\ell}=B_{i+1,\ell}$.
		Analogously to the case $\ell\in\{1,3\}$, with $x$ and $y$ interchanged,
		an edge gives $\phi(y)\le\phi_-(x)<\phi(x)$.
		For a nonedge, $\phi(x)<\phi_+(x)\le\phi(y)$.
		Hence $xy\in E(G)\Longleftrightarrow\phi(x)>\phi(y)$.
		
		For (ii), take $x_1\in A_d$ and $y_1\in Q_d$.
		Definition~$\ref{def:adjustment-values}$ gives $B^-_{d,2}=B_{d,1}=A_d$.
		Analogously to the case $\ell\in\{1,3\}$ in (i),
		an edge gives $\phi(x_1)\le\phi_-(y_1)<\phi(y_1)$.
		For a nonedge, $\phi(y_1)<\phi_+(y_1)\le\phi(x_1)$.
		Thus $x_1y_1\in E(G)\Longleftrightarrow\phi(x_1)<\phi(y_1)$.
		
		When $d\ge2$, take $x_2\in B_{d-1,3}$ and $y_2\in Q_{d-1}=B_{d-1,4}$.
		Since $B^-_{d-1,4}=B_{d-1,3}$, Definition~$\ref{def:adjustment-values}$
		and Lemma~$\ref{lem:adjustments}$\,(iv) give
		$\phi(x_2)\le\phi_-(y_2)<\phi(y_2)$ if $x_2y_2\in E(G)$ and
		$\phi(y_2)<\phi_+(y_2)\le\phi(x_2)$ if $x_2y_2\notin E(G)$.
		\nobreak\hfill$\Box$
		
		\subsection{Verification of the rank function}\label{subsec:rank-assignment}
		
		\noindent\textbf{Proof of Theorem~$\ref{the:sufficient}$.}
		Take two distinct vertices $w_1\in L_{i_1}$ and $w_2\in L_{i_2}$ with
		$0\le i_1\le i_2\le d$. Let $f_{\ell_1}(i_1)$ and
		$f_{\ell_2}(i_2)$ be their respective base values, where
		$w_j\in B_{i_j,\ell_j}$ for $w_j\ne u$, and $\ell_1=4$ when $w_1=u$.
		By Lemma~$\ref{lem:adjustments}$\,(ii),
		$|r(w_j)-f_{\ell_j}(i_j)|<\tfrac12$ for $w_j\ne u$, while
		$|r(u)-f_4(0)|=\tfrac12$.
		Since at most one of the two vertices is $u$,
		$|r(w_1)+r(w_2)-f_{\ell_1}(i_1)-f_{\ell_2}(i_2)|<1$.
		We consider the following three cases.
		
		\noindent\textbf{Case 1.} $i_2=i_1$.
		\\*\hspace*{\parindent}%
		Here $i_1\ge1$, since $L_0=\{u\}$. Each of $w_1$ and $w_2$ lies
		in the clique $Q_{i_1}$ or the stable set $A_{i_1}$.
		
		If $\{w_1,w_2\}\subseteq Q_{i_1}$, then $w_1w_2\in E(G)$.
		Since $Q_i=\emptyset$ for $0\le i\le d-2$, we have $i_1\in\{d-1,d\}$.
		By Lemma~$\ref{lem:adjustments}$\,(ii), $0<\phi(w_1)+\phi(w_2)<1$.
		For $i_1=d-1$, we have
		$r(w_1)+r(w_2)=\frac{4\theta_1}{3}+\phi(w_1)+\phi(w_2)
		\in(\frac{4\theta_1}{3},\frac{4\theta_1}{3}+1)\subsetneqq(\theta_1,\theta_2)$.
		For $i_1=d$, $f_2(d)=\theta_1$ gives
		$r(w_1)+r(w_2)=2\theta_1-\phi(w_1)-\phi(w_2)$.
		Since $\theta_1>3$ and $\theta_2>2\theta_1+2$, we obtain
		$\theta_1<2\theta_1-1<r(w_1)+r(w_2)<2\theta_1<\theta_2$.
		
		If $\{w_1,w_2\}\subseteq A_{i_1}$, then $w_1w_2\notin E(G)$.
		For $i_1=d$,
		$r(w_1)+r(w_2)=2(\theta_2-\theta_1)+\phi(w_1)+\phi(w_2)
		>2(\theta_2-\theta_1)>\theta_2$.
		For $i_1=d-1$, we have $\ell_1,\ell_2\le3$.
		Lemmas~$\ref{lem:base-sums}$\,(i)
		and~$\ref{lem:adjustments}$\,(ii) give
		$r(w_1)+r(w_2)<2f_3(d-1)+1=\frac{2\theta_1}{3}+1<\theta_1$.
		Suppose $i_1\le d-2$. If $d-i_1$ is odd, $f_4(i_1)$ is the largest
		base value by Lemma~$\ref{lem:base-sums}$\,(i). Equations~\eqref{eq:base-first} and~\eqref{100} give
		\[
		r(w_1)+r(w_2)<2f_4(i_1)+1
		\le-2(\theta_2-\theta_1)+\frac{4\theta_1}{3}+1<\theta_1.
		\]
		If $d-i_1$ is even, $f_4(i_1)$ is the smallest base value by
		Lemma~$\ref{lem:base-sums}$\,(i), so
		\[
		r(w_1)+r(w_2)>2f_4(i_1)-1
		\ge2(\theta_2-\theta_1)+\frac{2\theta_1}{3}-1>\theta_2.
		\]
		
		For one vertex in each part, by symmetry take
		$w_1\in A_{i_1}$ and $w_2\in Q_{i_1}$; then $i_1\in\{d-1,d\}$.
		If $i_1=d-1$ and $w_1\in B_{d-1,1}\cup B_{d-1,2}$,
		Lemma~$\ref{lem:part-adjacency}$\,(i) gives $w_1w_2\notin E(G)$, and
		$r(w_1)+r(w_2)<f_2(d-1)+f_4(d-1)+1
		=\frac{2\theta_1}{3}+1<\theta_1$.
		
		If $i_1=d-1$ and $w_1\in B_{d-1,3}$,
		Lemma~$\ref{lem:base-sums}$\,(iv) gives
		$r(w_1)+r(w_2)=\theta_1-\phi(w_1)+\phi(w_2)$.
		For an edge, Lemma~$\ref{lem:adjustment-adjacency}$\,(ii) gives
		$\phi(w_1)<\phi(w_2)$; together with Lemma~$\ref{lem:adjustments}$\,(ii),
		$\theta_1<r(w_1)+r(w_2)<\theta_1+\tfrac12<\theta_2$.
		For a nonedge, Lemma~$\ref{lem:adjustment-adjacency}$\,(ii) gives
		$\phi(w_2)<\phi(w_1)$, so $r(w_1)+r(w_2)<\theta_1$.
		For $i_1=d$, $r(w_1)+r(w_2)=\theta_2+\phi(w_1)-\phi(w_2)$.
		Analogously to the $B_{d-1,3},Q_{d-1}$ case, using
		Lemma~$\ref{lem:adjustment-adjacency}$\,(ii) for $A_d,Q_d$,
		an edge gives $\theta_1<\theta_2-\tfrac12<r(w_1)+r(w_2)<\theta_2$,
		and a nonedge gives $r(w_1)+r(w_2)>\theta_2$.
		
		\noindent\textbf{Case 2.} $i_2=i_1+1$.
		\\*\hspace*{\parindent}%
		If $i_1=0$, then $w_1=u$ and $w_1w_2\in E(G)$.
		For $\ell_2<4=\ell_1$, Lemmas~$\ref{lem:base-sums}$\,(i)--(iii)
		and~$\ref{lem:adjustments}$\,(ii) give
		$\theta_1<f_4(0)+f_{\ell_2}(1)-1
		<r(u)+r(w_2)<f_4(0)+f_{\ell_2}(1)+1<\theta_2$.
		If $\ell_2=4$ and $d$ is odd, Lemmas~$\ref{lem:base-sums}$\,(iii)
		and~$\ref{lem:adjustments}$\,(ii) give
		\[
		r(u)+r(w_2)=\theta_1+\tfrac12-\phi(w_2)
		\in(\theta_1,\theta_1+\tfrac12)\subsetneqq(\theta_1,\theta_2).
		\]
		If $\ell_2=4$ and $d$ is even, Lemmas~$\ref{lem:base-sums}$\,(iii)
		and~$\ref{lem:adjustments}$\,(ii) give
		\[
		r(u)+r(w_2)=\theta_2-\tfrac12+\phi(w_2)
		\in(\theta_2-\tfrac12,\theta_2)\subsetneqq(\theta_1,\theta_2).
		\]
		
		Suppose $i_1\ge1$.
		If $\ell_2<\ell_1$, Lemma~$\ref{lem:part-adjacency}$\,(iii) gives
		$w_1w_2\in E(G)$.
		As in the root case with $\ell_2<4$, Lemma~$\ref{lem:base-sums}$\,(i)--(iii) gives
		$\theta_1+1<f_{\ell_1}(i_1)+f_{\ell_2}(i_1+1)<\theta_2-1$.
		Lemma~$\ref{lem:adjustments}$\,(ii) gives
		$\theta_1<r(w_1)+r(w_2)<\theta_2$.
		If $\ell_2>\ell_1$, Lemma~$\ref{lem:part-adjacency}$\,(iii) gives
		$w_1w_2\notin E(G)$.
		If $d-i_1$ is odd, Lemma~$\ref{lem:base-sums}$\,(i)--(iii) gives
		\[
		r(w_1)+r(w_2)<f_{\ell_1}(i_1)+f_{\ell_2}(i_1+1)+1
		=\theta_1+f_{\ell_1}(i_1)-f_{\ell_2}(i_1)+1<\theta_1.
		\]
		If $d-i_1$ is even, Lemma~$\ref{lem:base-sums}$\,(i)--(iii) gives
		\[
		r(w_1)+r(w_2)>f_{\ell_1}(i_1)+f_{\ell_2}(i_1+1)-1
		=\theta_2+f_{\ell_1}(i_1)-f_{\ell_2}(i_1)-1>\theta_2.
		\]
		
		If $\ell_1=\ell_2=\ell\in\{1,3\}$ and $d-i_1$ is odd,
		Lemma~$\ref{lem:base-sums}$\,(iii) gives
		$r(w_1)+r(w_2)=\theta_1-\phi(w_1)+\phi(w_2)$.
		For an edge, Lemma~$\ref{lem:adjustment-adjacency}$\,(i) gives
		$\phi(w_1)<\phi(w_2)$, so
		$\theta_1<r(w_1)+r(w_2)<\theta_1+\tfrac12<\theta_2$.
		For a nonedge, Lemma~$\ref{lem:adjustment-adjacency}$\,(i) gives $\phi(w_2)<\phi(w_1)$,
		so $r(w_1)+r(w_2)<\theta_1$.
		
		If $\ell_1=\ell_2=\ell\in\{1,3\}$ and $d-i_1$ is even,
		$r(w_1)+r(w_2)=\theta_2+\phi(w_1)-\phi(w_2)$.
		By Lemma~$\ref{lem:adjustment-adjacency}$\,(i), an edge gives $\theta_1<\theta_2-\tfrac12<r(w_1)+r(w_2)<\theta_2$,
		and a nonedge gives $r(w_1)+r(w_2)>\theta_2$.
		
		If $\ell_1=\ell_2=\ell\in\{2,4\}$ and $d-i_1$ is odd,
		$r(w_1)+r(w_2)=\theta_1+\phi(w_1)-\phi(w_2)$.
		By Lemma~$\ref{lem:adjustment-adjacency}$\,(i), an edge gives
		$\phi(w_1)>\phi(w_2)$, so $\theta_1<r(w_1)+r(w_2)<\theta_1+\tfrac12<\theta_2$,
		and a nonedge gives $r(w_1)+r(w_2)<\theta_1$.
		
		If $\ell_1=\ell_2=\ell\in\{2,4\}$ and $d-i_1$ is even,
		$r(w_1)+r(w_2)=\theta_2-\phi(w_1)+\phi(w_2)$.
		By Lemma~$\ref{lem:adjustment-adjacency}$\,(i), an edge gives
		$\phi(w_1)>\phi(w_2)$, so $\theta_1<\theta_2-\tfrac12<r(w_1)+r(w_2)<\theta_2$,
		and a nonedge gives $r(w_1)+r(w_2)>\theta_2$.
		
		\noindent\textbf{Case 3.} $i_2\ge i_1+2$.
		\\*\hspace*{\parindent}%
		Since $d_G(u,w_2)-d_G(u,w_1)=i_2-i_1\ge2$, we have $w_1w_2\notin E(G)$.
		The assumption $\theta_2>2\theta_1+2$ gives $\theta_2-\theta_1>\theta_1+2$.
		Recall that $\ell_2\in\{1,2\}$ when $i_2=d$.
		
		If $d-i_1=2k+1$, then $k\ge1$ and $d-i_2\le2k-1$.
		When $d-i_2$ is even, $d-i_2\le2k-2$, so Lemma~$\ref{lem:base-sums}$\,(i)
		and equation~\eqref{eq:base-first} give
		$f_{\ell_2}(i_2)\le f_1(i_2)\le k(\theta_2-\theta_1)$.
		When $d-i_2$ is odd, $f_{\ell_2}(i_2)\le f_4(i_2)$ by Lemma~$\ref{lem:base-sums}$\,(i),
		and equations~\eqref{eq:base-first} and~\eqref{100} give $f_4(i_2)\le\tfrac{2\theta_1}{3}$.
		Thus
		\[
		f_{\ell_2}(i_2)\le\max\{k(\theta_2-\theta_1),\tfrac{2\theta_1}{3}\}=k(\theta_2-\theta_1).
		\]
		Since $d-i_1$ is odd, Lemma~$\ref{lem:base-sums}$\,(i) gives
		$f_{\ell_1}(i_1)\le f_4(i_1)$, so
		\[
		f_{\ell_1}(i_1)+f_{\ell_2}(i_2)
		\le f_4(i_1)+k(\theta_2-\theta_1)
		=-k(\theta_2-\theta_1)+\tfrac{2\theta_1}{3}+k(\theta_2-\theta_1)=\tfrac{2\theta_1}{3}.
		\]
		Lemma~$\ref{lem:adjustments}$\,(ii) now gives
		\[
		r(w_1)+r(w_2)<f_{\ell_1}(i_1)+f_{\ell_2}(i_2)+1
		\le\frac{2\theta_1}{3}+1<\theta_1.
		\]
		
		For $d-i_1=2k$, $k\ge1$, we have $d-i_2\le2k-2$.
		If $d-i_2$ is even, then $f_{\ell_2}(i_2)\ge\theta_1$ since
		$\ell_2\in\{1,2\}$ when $i_2=d$, and
		$f_{\ell_2}(i_2)\ge f_4(i_2)\ge(\theta_2-\theta_1)+\tfrac{\theta_1}{3}>\theta_1$ otherwise.
		If $d-i_2$ is odd, equations~\eqref{eq:base-first} and~\eqref{100} give
		$f_{\ell_2}(i_2)\ge f_1(i_2)\ge\theta_1-(k-1)(\theta_2-\theta_1)$.
		Thus
		\[
		f_{\ell_2}(i_2)\ge\min\{\theta_1,\theta_1-(k-1)(\theta_2-\theta_1)\}
		=\theta_1-(k-1)(\theta_2-\theta_1).
		\]
		Since $d-i_1$ is even, Lemma~$\ref{lem:base-sums}$\,(i) gives
		$f_{\ell_1}(i_1)\ge f_4(i_1)$, so
		\[
		f_{\ell_1}(i_1)+f_{\ell_2}(i_2)
		\ge f_4(i_1)+\theta_1-(k-1)(\theta_2-\theta_1)
		=k(\theta_2-\theta_1)+\tfrac{\theta_1}{3}+\theta_1-(k-1)(\theta_2-\theta_1)
		=\theta_2+\tfrac{\theta_1}{3}.
		\]
		Lemma~$\ref{lem:adjustments}$\,(ii) therefore gives
		\[
		r(w_1)+r(w_2)>f_{\ell_1}(i_1)+f_{\ell_2}(i_2)-1
		\ge\theta_2+\frac{\theta_1}{3}-1>\theta_2.
		\]
		
		Thus the rank sum of every edge lies in $(\theta_1,\theta_2)$, and
		the rank sum of every nonedge lies outside $[\theta_1,\theta_2]$.
		\nobreak\hfill$\Box$
		
		\Needspace{4\baselineskip}
		\begin{samepage}
			\begin{corollary}\label{the:main}
				A connected non-bipartite graph is $2$-threshold if and only
				if it admits a BFS root and level decompositions satisfying
				conditions $(1)$--$(5)$ of Theorem~$\ref{the:necessary}$.
			\end{corollary}
		\end{samepage}

		\section{Bipartite graphs and disconnected graphs}
		\subsection{The bipartite case}\label{subsec:bipartite}
		Jamison and Sprague \cite{Jamison1} proved that bipartite
		$2$-threshold graphs are precisely bipartite permutation graphs.
		The following theorem describes this class using BFS levels.
		
		\begin{theorem}\label{the:bipartite}
			Let $G$ be a connected bipartite graph with at least two vertices.
			Then $G$ is $2$-threshold if and only if it has a BFS root $u$
			with levels $L_0,\ldots,L_d$ satisfying the following conditions.
			
			\textup{(i)} For each $0\le i<d$, the family
			$\{N_{L_{i+1}}(v):v\in L_i\}$ is nested.
			
			\textup{(ii)} For $1\le i<d$, take $x\in D_{i,t}$ and
			$y\in D_{i,t'}$ with $0\le t<t'\le s_i$.
			Then $N_{L_{i-1}}(y)\subseteq N_{L_{i-1}}(x)$.
			Here the classes $D_{i,t}$ are defined as in
			Definition~$\ref{def:next-classes}$ using the nesting in \textup{(i)}.
		\end{theorem}
		
		\noindent\textbf{Proof.}
		Suppose first that $G$ is $2$-threshold.
		Choose a strict representation by Lemma~$\ref{lem:strict}$ and
		let $u$ be a vertex of maximum rank.
		Order odd levels by increasing rank and even levels by decreasing rank.
		Corollary~$\ref{cor:next-nesting}$\,(ii) gives condition \textup{(i)}.
		For \textup{(ii)}, take $x_1\in D_{i,t}$ and $y_1\in D_{i,t'}$ with
		$1\le i<d$ and $0\le t<t'\le s_i$.
		Definition~$\ref{def:next-classes}$ gives
		$N_{L_{i+1}}(x_1)\subsetneqq N_{L_{i+1}}(y_1)$, so $x_1$ precedes $y_1$.
		Corollary~$\ref{cor:next-nesting}$\,(ii) then gives
		$N_{L_{i-1}}(y_1)\subseteq N_{L_{i-1}}(x_1)$, proving \textup{(ii)}.
		
		Conversely, take a BFS root $u$ satisfying \textup{(i)} and \textup{(ii)}.
		Each $L_i$ lies on one side of the bipartition, so $A_i=L_i$ and
		$Q_i=\emptyset$ give an admissible decomposition for $0\le i\le d$.
		In Theorem~$\ref{the:necessary}$, the sets $Q_i$ in $(4)$ and the sets
		$R_t^{(i)}\subseteq Q_i$ in $(5)$ are empty. To apply
		Theorem~$\ref{the:sufficient}$, it remains to check $(1)$--$(3)$.
		
		By \textup{(i)}, the neighborhoods from each level into the next are
		nested, as required in $(1)$. For $(2)$, the required inclusion holds since
		$Q_{d-1}\cup N_{A_{d-1}}(Q_{d-1})=\emptyset\subseteq D_{d-1,s_{d-1}}$.
		For $1\le i\le d-2$, \textup{(ii)} gives $(3)$.
		When $i=d-1\ge1$, all vertices of $D_{i,s_i}$ have empty neighborhood
		in $L_i$. Definition~$\ref{def:stable-classes}$\,(i) gives $m_i=s_i$ and
		$S_{i,t}=D_{i,t}$ for $0\le t\le s_i$. Thus \textup{(ii)} gives
		$N_{L_{i-1}}(y)\subseteq N_{L_{i-1}}(x)$ for $x\in S_{i,t}$,
		$y\in S_{i,t'}$ with $t<t'$. At $i=d$, Definition~$\ref{def:stable-classes}$\,(ii) gives
		$S_{d,0}=L_d$ and $m_d=0$, so no pair $(t,t')$ satisfies
		$0\le t<t'\le m_d$ in $(3)$.
		Theorem~$\ref{the:sufficient}$ now gives a $2$-threshold representation
		of $G$. $\Box$
		
		\Needspace{6\baselineskip}
		\subsection{The disconnected case}
		For arbitrary graphs, we use the component criterion of Xiao
		and Nagamochi \cite{Xiao}.
		
		\begin{theorem}[{\cite{Xiao}}]\label{the:components}
			A graph is a $2$-threshold graph if and only if it has at most one
			non-bipartite component and all of its components are $2$-threshold
			graphs.
		\end{theorem}
		
		\Needspace{3\baselineskip}
		A single vertex is $2$-threshold. Applying Corollary~$\ref{the:main}$ and
		Theorem~$\ref{the:bipartite}$ to the components in
		Theorem~$\ref{the:components}$ gives the following characterization.
		
		\begin{corollary}
			A graph $G$ is a $2$-threshold graph if and only if\nopagebreak
			
			\textup{(i)} at most one component is non-bipartite;
			
			\textup{(ii)} every bipartite component with at least two vertices has
			a BFS root satisfying conditions \textup{(i)} and \textup{(ii)} of
			Theorem~$\ref{the:bipartite}$;
			
			\textup{(iii)} the non-bipartite component, if present, has a BFS root
			and admissible decompositions satisfying conditions $(1)$--$(5)$
			of Theorem~$\ref{the:necessary}$.
		\end{corollary}

		\section{Concluding Remarks}\label{sec:conclusion}
		We characterize connected $2$-threshold graphs by conditions on BFS levels
		and construct a rank representation whenever these conditions hold.
		Using Xiao and Nagamochi's theorem on connected components, we also characterize
		disconnected $2$-threshold graphs.
		
		Future work could extend the BFS characterization to connected
		$k$-threshold graphs for $k\ge3$.
		
		\section*{Declaration of AI use}
		The authors used ChatGPT (OpenAI) to assist with language editing and manuscript revision.
		The authors reviewed the revisions and take full responsibility for the content of the manuscript.

	\end{sloppypar}
\end{document}